\documentclass[reqno]{amsart}
\usepackage[centertags]{amsmath}
\usepackage{graphicx}
\usepackage{amsfonts}
\usepackage{amssymb}
\usepackage{amsthm}
\usepackage[top=1in, bottom=1in, left=1in, right=1in]{geometry}

\usepackage[pdfborder={0 0 0},
   pdftitle={SUBDIVISION RULES AND THE EIGHT MODEL GEOMETRIES},
  pdfauthor={BRIAN RUSHTON}, pdftex, bookmarks=true, bookmarksnumbered = true]{hyperref}

\theoremstyle{definition}
\newtheorem{thm}{Theorem}
\theoremstyle{definition}

\theoremstyle{definition}
\newtheorem*{defi}{Definition}
\theoremstyle{remark}

\theoremstyle{definition}

\newcommand{\R}{\mathbb{R}}

\newcommand{\parlengths}{\setlength{\parindent}{0pt}}
\begin{document}
\date{\today}

\pdfbookmark[1]{SUBDIVISION RULES AND THE EIGHT MODEL GEOMETRIES}{user-title-page}

\title{Subdivision rules and the eight model geometries}

\author{Brian Rushton}
\address{Department of Mathematics, Brigham Young University, Provo, UT 84602, USA}
\email{lindianr@gmail.com}

\begin{abstract}
Cannon and Swenson have shown that each hyperbolic 3-manifold group has a natural subdivision rule on the space at infinity \cite{hyperbolic}, and that this subdivision rule captures the action of the group on the sphere. Explicit subdivision rules have also been found for some close finite-volume hyperbolic manifolds, as well as a few non-hyperbolic knot complements \cite{linksubs}, \cite{myself2}. We extend these results by finding explicit finite subdivision rules for closed manifolds of the $\mathbb{E}^3$, $\mathbb{H}^2 \times \mathbb{R}$, $\mathbb{S}^2 \times \mathbb{R}$, $\mathbb{S}^3$, and $\widetilde{SL_2(\mathbb{R})}$ manifolds by means of model manifolds. Because all manifolds in these geometries are the same up to finite covers, the subdivision rules for these model manifolds will be very similar to subdivision rules for all other manifolds in their respective geometries. We also discuss the existence of subdivision rules for Nil and Sol geometries. We use Ken Stephenson's Circlepack \cite{Circlepak} to visualize the subdivision rules and the resulting space at infinity.
\end{abstract}

\maketitle\parlengths

\begin{figure}
\begin{center}
\scalebox{.65}{\includegraphics{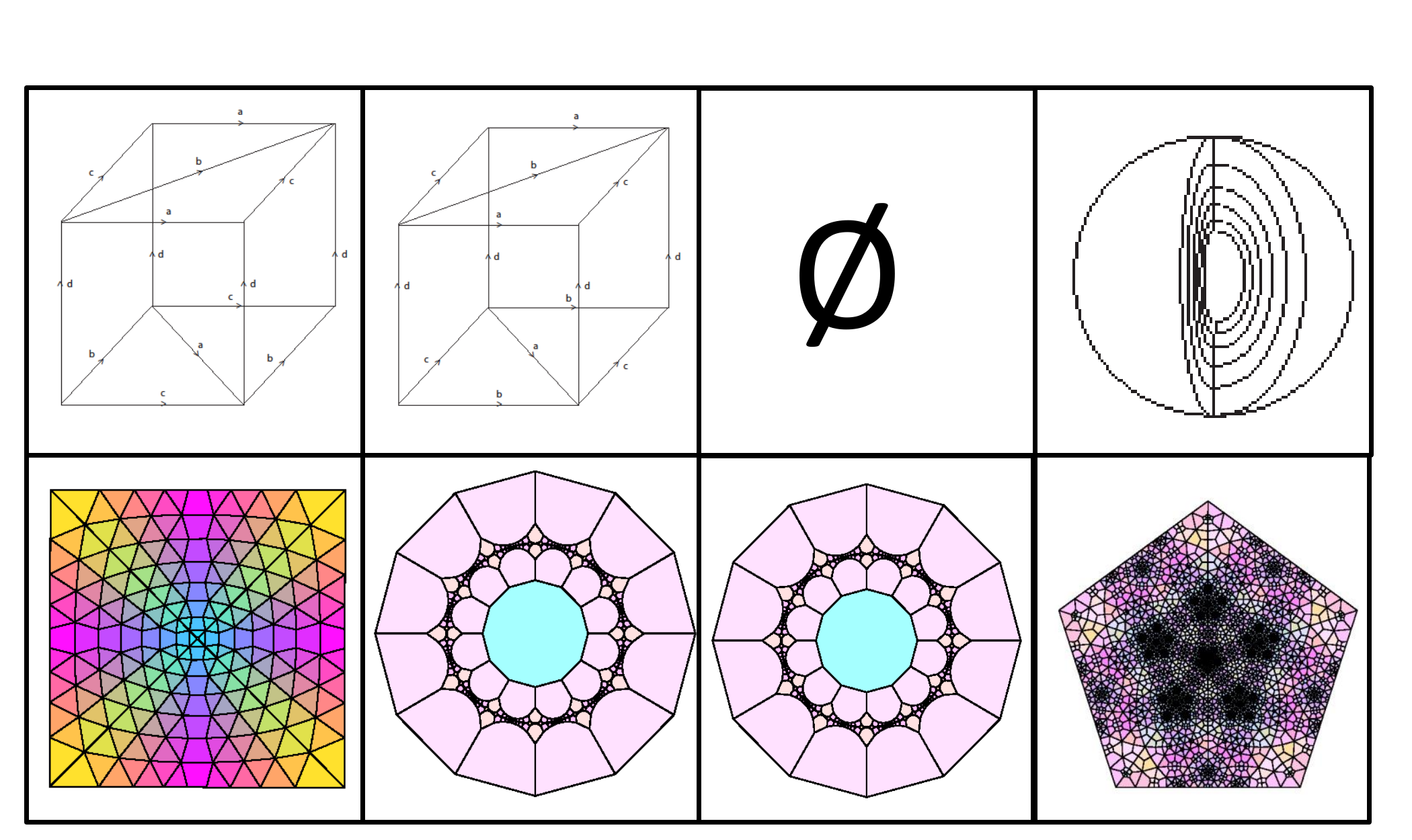}} \caption{Representations of the eight geometries. The first two figures represent the most basic Nil and Sol manifolds, which do not have finite subdivision rules. The other six are subdivision complexes for manifolds of the $\mathbb{S}^3$ and $\mathbb{S}^2\times \mathbb{R}$, $\mathbb{E}^3$, $\mathbb{H}^2 \times \mathbb{R}$, $\widetilde{SL_2(\mathbb{R})}$, and $\mathbb{H}^3$ geometries. See section \ref{summary}.}
\label{AllSubs}
\end{center}
\end{figure}

\section{Background}
This paper assumes basic knowledge of algebraic topology (fundamental group, covering spaces, etc.) and Thurston's eight geometries \cite{Thurston}. Peter Scott's exposition of the eight geometries \cite{Scott} is especially helpful, and will be our main reference for properties of the geometries.
\section{Introduction}

One of the most useful characteristics of hyperbolic 3-manifold groups (or of Gromov hyperbolic groups in general) is their space at infinity. One construction related to the space at infinity is a subdivision rule. A subdivision rule is a recursive way of dividing polygons (or, more generally, compact sets that cover a space) into smaller polygons (or compact sets) according to finitely many rules or tile types; for instance, barycentric subdivision is a subdivision rule with a single tile type (a triangle), where each tile is subdivided into six smaller tiles and each edge is divided in two.

Cannon and Swenson have shown \cite{hyperbolic} that every hyperbolic group with a 2-sphere at infinity has a subdivision rule which acts on the sphere at infinity. Cannon, Floyd, and Parry have studied such subdivision rules extensively in an attempt to prove a Cannon's conjecture, which states that every Gromov hyperbolic group with a 2-sphere at infinity acts cocompactly and properly discontinuously by isometries on hyperbolic 3-space. Cannon's conjecture can be proven if it can be shown that the subdivision rule of the sphere at infinity is conformal in a certain sense \cite{Combinatorial}.

\begin{figure}
\begin{center}
\scalebox{.65}{\includegraphics{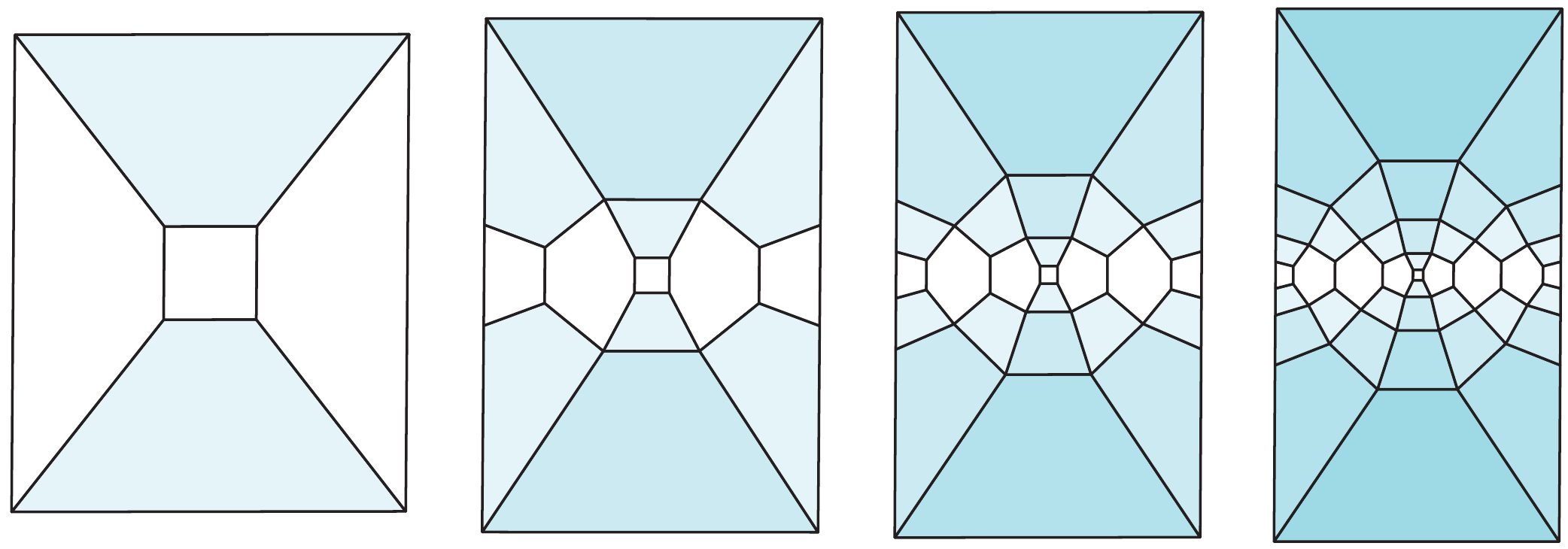}} \caption{Several subdivisions of a subdivision complex associated to the Hopf link complement.}
\label{HopfBound}
\end{center}
\end{figure}

\begin{figure}
\begin{center}
\scalebox{.65}{\includegraphics{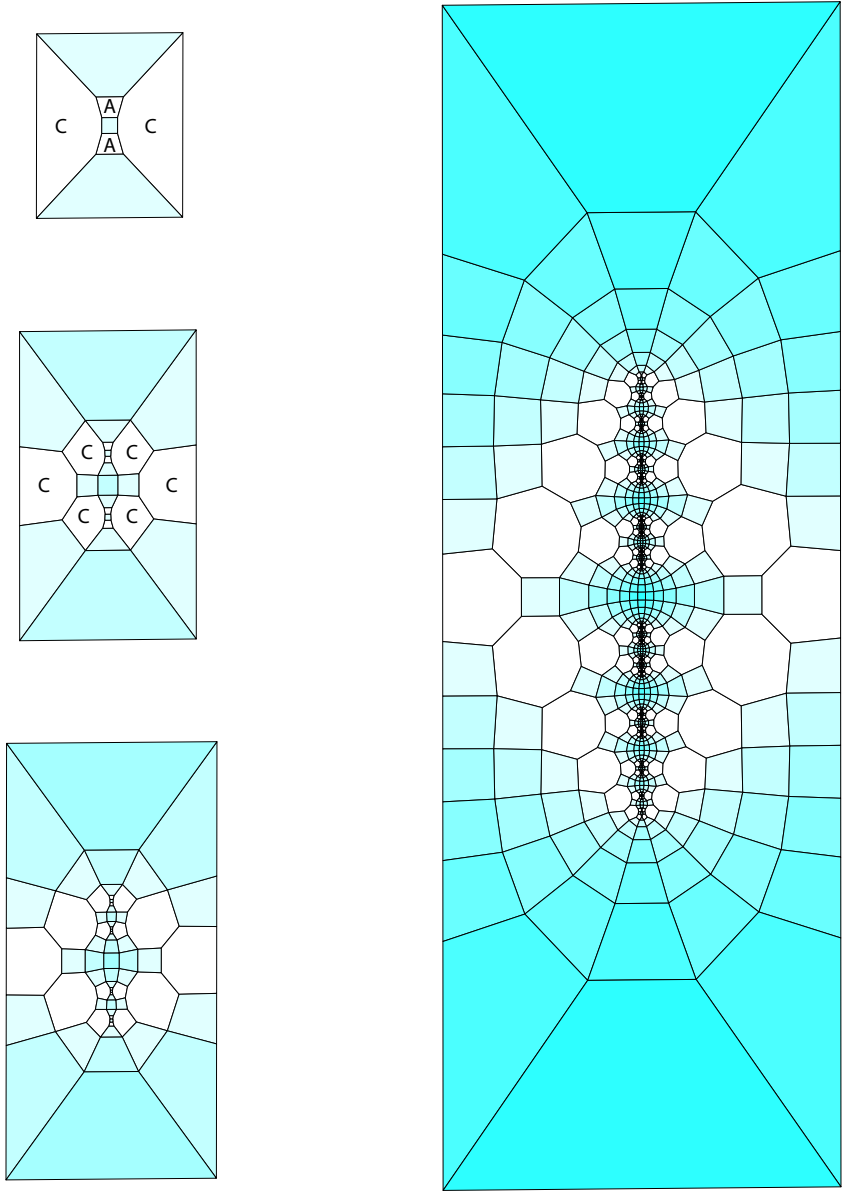}} \caption{Several subdivisions of a subdivision complex associated to the trefoil knot complement.}
\label{trefbound}
\end{center}
\end{figure}

\begin{figure}
\begin{center}
\scalebox{.65}{\includegraphics{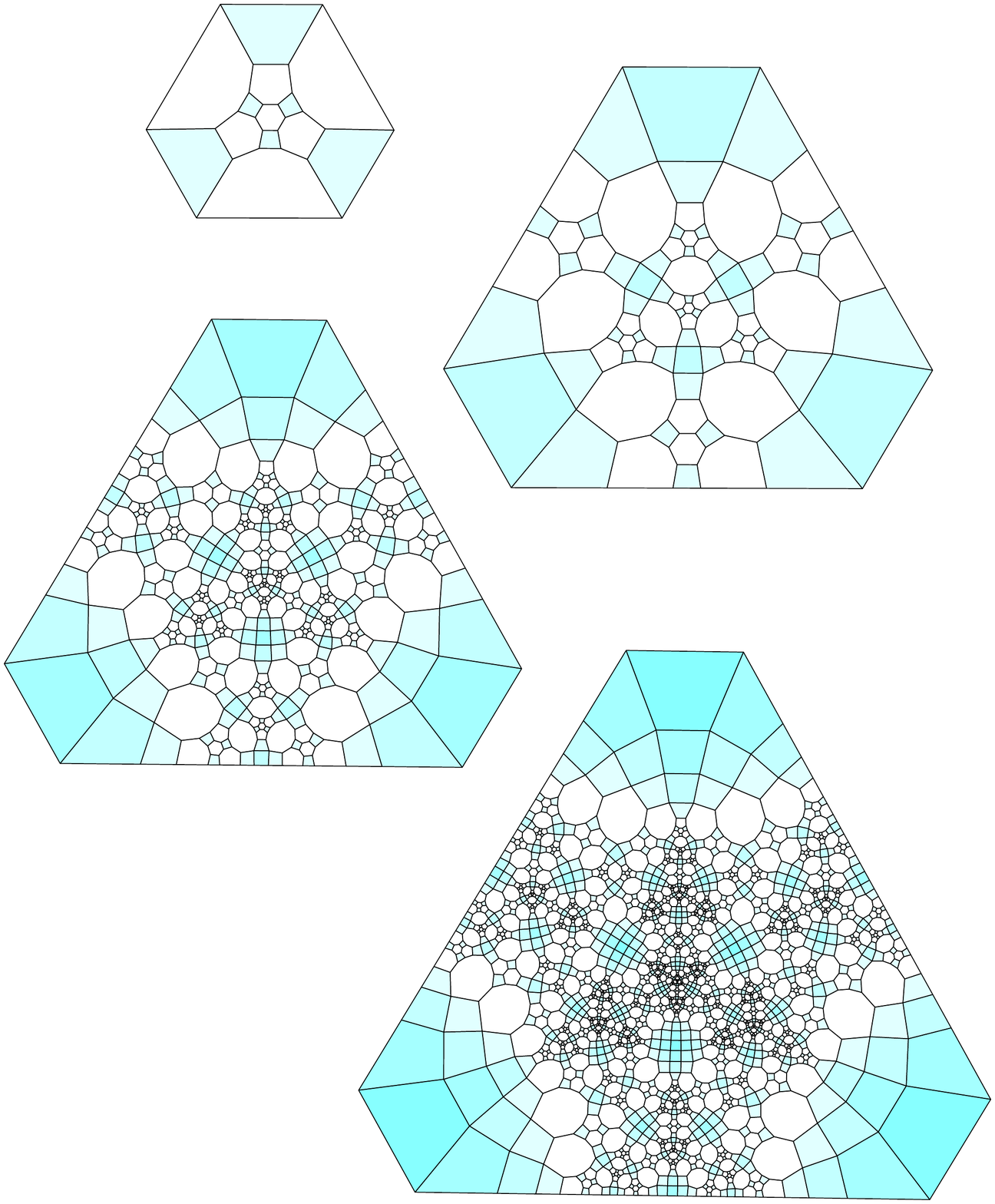}} \caption{Several subdivisions of a subdivision complex associated to the Borromean rings complement.}
\label{Bringbound}
\end{center}
\end{figure}

Thus, subdivision rules are intimately connected with hyperbolic groups. It was surprising, then, when Cannon discovered a finite subdivision rule associated to the trefoil knot, which is not a hyperbolic knot. His example, though unpublished, was later expanded on by the author in \cite{linksubs}, in which we find finite subdivision rules (defined in the next section) for all non-split, prime, alternating link complements, including the non-hyperbolic 2-braid links. The geometry of the link complements was apparent in the subdivision rules; the complement of the Hopf link (which is a Euclidean manifold) had linear growth in number of cells, while all other 2-braid link complements (which belong to the $\mathbb{H^2}\times\mathbb{R}$ geometry) grew exponentially in one direction and linearly in the other. The other prime alternating link complements (which are finite volume hyperbolic manifolds) divided exponentially in every direction. These three examples are shown in Figures \ref{HopfBound}, \ref{trefbound}, and \ref{Bringbound}.

This connection between the geometry of a 3-manifold and its subdivision rule led us to consider subdivision rules for closed 3-manifolds of other geometries. Surprisingly, we have the following result:

\begin{thm}\label{BigTheorem}
There exists a finite subdivision rule that constructs the universal cover for at least one manifold in each of the $\mathbb{S}^3$ and $\mathbb{S}^2\times \mathbb{R}$, $\mathbb{E}^3$, $\mathbb{H}^2 \times \mathbb{R}$, $\widetilde{SL_2(\mathbb{R})}$, and $\mathbb{H}^3$ geometries, and no finite subdivision rule can construct the universal cover of a Sol manifold.
\end{thm}

(Note that Nil manifolds are not discussed in the theorem. See Section \ref{NilSection} for a discussion of this geometry.)

The next seven sections constitute the proof of this theorem.

Just like alternating link complements, the character of the subdivision rule depends on the geometry of the closed 3-manifold (this can be seen in Figure \ref{AllSubs}. The construction of these subdivision rules will be the focus of the main portion of the paper, followed by a discussion of the other geometries, Nil and Sol. At the beginning of each appropriate section, we show a circle-packed image of several stages of subdivision created with Ken Stephenson's Circlepack \cite{Circlepak}.

The subdivision rules we obtain are not at all unique, and are highly dependent on the choice of generating set. However, subdivision rules for all closed manifolds in a fixed geometry share certain characteristics. This topic requires a different tools to study than those used in this paper, and will be the subject of future work.

\section{Finite subdivision rules and finite replacement rules}

While our intuitive idea of a subdivision rule is helpful, we need a concrete definition of subdivision rule. We first recall Cannon, Floyd and Parry's definition of a finite subdivision rule, taken from \cite{subdivision}.

\begin{defi} A \textbf{finite subdivision rule} $R$ consists of the following.
\begin{enumerate}
\item A finite 2-dimensional CW complex $S_R$, called the \textbf{subdivision complex}, with a fixed cell structure such that $S_R$ is the union of its closed 2-cells. We assume that for each closed 2-cell $\tilde{s}$ of $S_R$ there is a CW structure $s$ on a closed 2-disk such that $s$ has at least three vertices, the vertices and edges of $s$ are contained in $\partial s$, and the characteristic map $\psi_s:s\rightarrow S_R$ which maps onto $\tilde{s}$ restricts to a homeomorphism onto each open cell.
\item A finite two dimensional CW complex $R(S_R)$, which is a subdivision of $S_R$.
\item A continuous cellular map $\phi_R:R(S_R)\rightarrow S_R$ called the \textbf{subdivision map}, whose restriction to every open cell is a homeomorphism.
\end{enumerate}
\end{defi}

Each CW complex $s$ in the definition above (with its given characteristic map $\psi_s$) is called a \textbf{tile type}.

As the final part of the definition, they show how finite subdivision rules can act on surfaces (and 2-complexes in general). An $R$-complex for a subdivision rule $R$ is a 2-dimensional CW complex $X$ which is the union of its closed 2-cells, together with a continuous cellular map $f:X\rightarrow S_R$ whose restriction to each open cell is a homeomorphism. We can subdivide $X$ into a complex $R(X)$ by requiring that the induced map $f:R(X)\rightarrow R(S_R)$ restricts to a homeomorphism onto each open cell. $R(X)$ is again an $R$-complex with map $\phi_R \circ f:R(X)\rightarrow S_R$. By repeating this process, we obtain a sequence of subdivided $R$-complexes $R^n(X)$ with maps $\phi_R^n\circ f:R^n(X)\rightarrow S_R$. All of the preceding definitions were adapted from \cite{subdivision}, which contains several examples. While in theory, a subdivision rule is represented by a CW-complex, most rules in practice are described by diagrams (such as Figure \ref{TorusSubs}).

To create a subdivision rule for a manifold, we first need a sphere for the subdivision rule to act on. Let $B(0)$ be a polyhedral fundamental domain for a manifold, and let $S(0)$ be its boundary, with the cell structure it inherits from $B(0)$.  Now, let $B(1)$ be formed from $B(0)$ by attaching polyhedra to all its open faces, and let $S(1)$ be its boundary, and so on.  For many manifolds and many choices of fundamental domain, $S(n)$ will always be a 2-sphere,and this process defines a sequence of tilings of the sphere.

We'd like to find a recursive way of describing this space at infinity, to lead us towards a subdivision rule.  In particular, we're looking for a replacement rule.

\begin{defi}\label{ReplaceDef} A \textbf{finite replacement rule} is a finite subdivision rule coupled with a \textbf{local combination rule}. In a local combination rule, a labeled tiling $T(n)$ of a surface is made more coarse by removing some edges and vertices in a local manner. More specifically, the set of tiles in $T(n)$ is partitioned into a number of closed sets $F_i$ with disjoint interior which are topological closed disks, each the union of closed labeled tiles. Each $F_i$ is replaced with a single closed labeled tile whose boundary is that of $F_i$. The sets $F_i$ may be of size one (meaning that the tile is not combined with any other tile). The partition and the replacement rule are completely determined by the labels of the tiles and of their neighbors of distance $<N$ for a fixed constant $N$ (here, the distance between two tiles $T, T'$ is given by the shortest length of a chain of closed tiles $T_1=T, T_2,...,T_k=T'$ such that $\cup T_i$ is connected). Thus, if there are two tiles in the surface with sufficiently large neighborhoods $A, A'$ that are isomorphic (meaning they are cellularly isomorphic and the isomorphism preserves labels), there are two smaller closed neighborhoods $B \subseteq A, B' \subseteq A'$ that are isomorphic in the post-combination tiling.

A finite replacement rule, then, acts on a complex by alternating subdivisions according to a finite subdivision rule and combinations according to a local combination rule.
\end{defi}

To turn a replacement rule into a subdivision rule, we need to eliminate the combination step above.

In the sections that follow, we will examine each of Thurston's eight geometries, one in each section. We pick a representative manifold, find a subdivision rule where possible, and display circle packings of the subdivision rules using Ken Stephenson's Circlepack \cite{Circlepak}. For five of the geometries, the representative manifold will be essentially the only manifold. In the $\mathbb{E}^3$, $\mathbb{H}^2 \times \mathbb{R}$, and $\mathbb{S}^2 \times \mathbb{R}$ geometries, every manifold is, up to a finite cover, a trivial circle bundle over a Euclidean, hyperbolic, or spherical 2-manifold, respectively (see \cite{Scott}). We will use these trivial bundles as our model manifolds for these geometries.

All $\widetilde{SL_2(\mathbb{R}}$ manifolds are, up to finite covers, unit tangent bundles over hyperbolic surfaces, which all cover the unit tangent bundle of the smallest hyperbolic surface. All $\mathbb{S}^3$ manifolds are finitely covered by $\mathbb{S}^3$ \cite{Scott}. Thus, these two geometries can also each be represented by a single closed 3-manifold, which we examine. Similar statements hold for Nil and Sol, but we will see that we have no subdivision rules for these geometries.

The reason we are only interested in manifolds up to finite covers is that, in geometric group theory, finite index subgroups are almost (or `virtually') the same as the original group. And, as experimentation will show, similar groups have similar subdivision rules (the 2-braid knots with more than 2 twists in \cite{linksubs} are a good family of examples).

As mentioned earlier, these subdivision rules are not unique. However, there are some characteristics that subdivision rules in a particular geometry must share; for instance, all of the virtually abelian groups must have polynomial growth.

\section{\texorpdfstring{$E^3$ geometry: $S^1\times S^1 \times S^1$}{E\textthreesuperior\ geometry:\ S\textonesuperior\ x S\textonesuperior\ x S\textonesuperior}}\label{ThreeTorus}

\begin{figure}
\begin{center}
\scalebox{.8}{\includegraphics{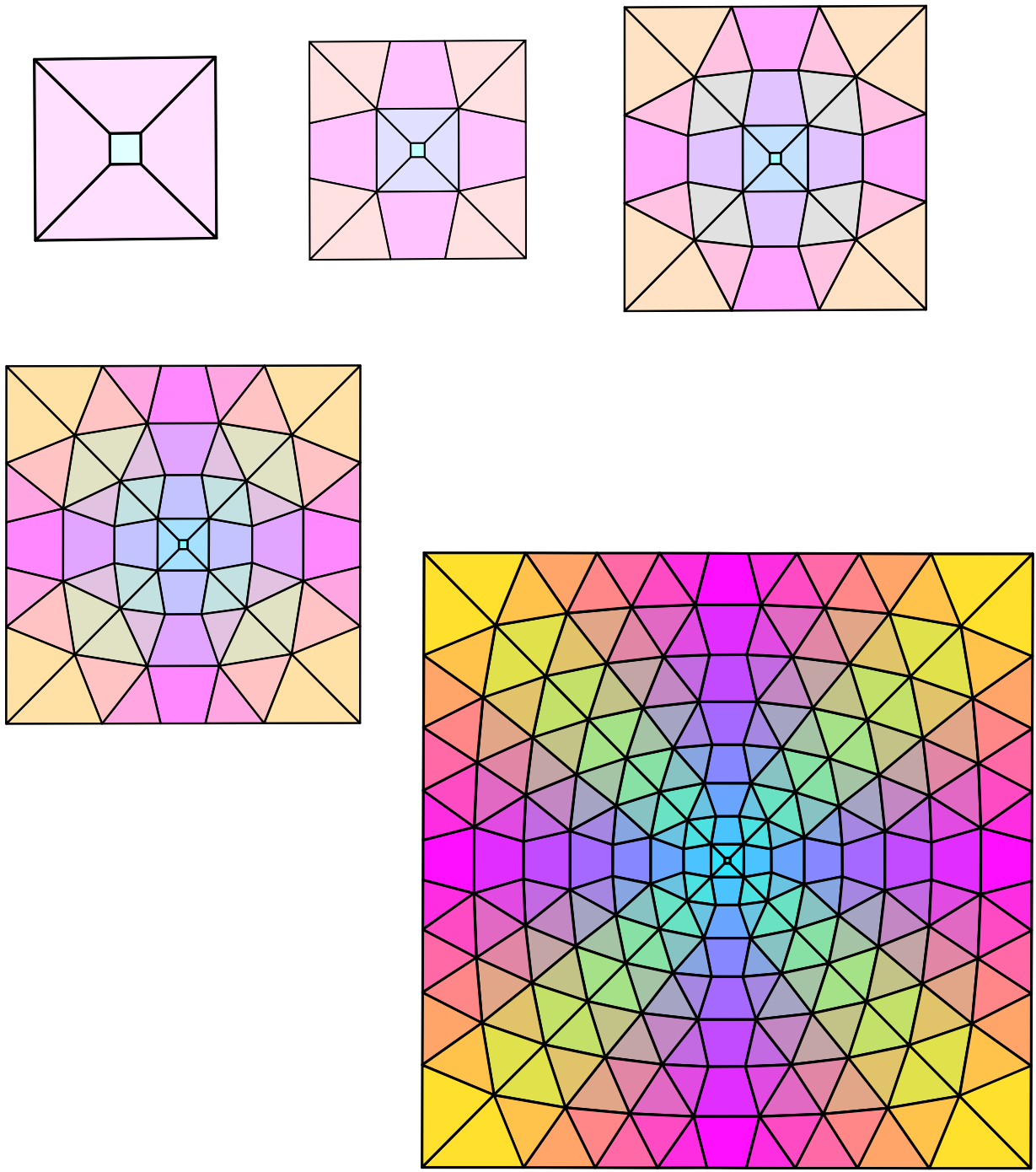}} \caption{The subdivision complexes $R(X), R^2 (X), R^3(X), R^4(X),$ and $R^8(X)$, where $X$ is a face of a fundamental domain of $S^1\times S^1 \times S^1$.}
\label{CircleTorus}
\end{center}
\end{figure}

\begin{figure}
\begin{center}
\scalebox{.8}{\includegraphics{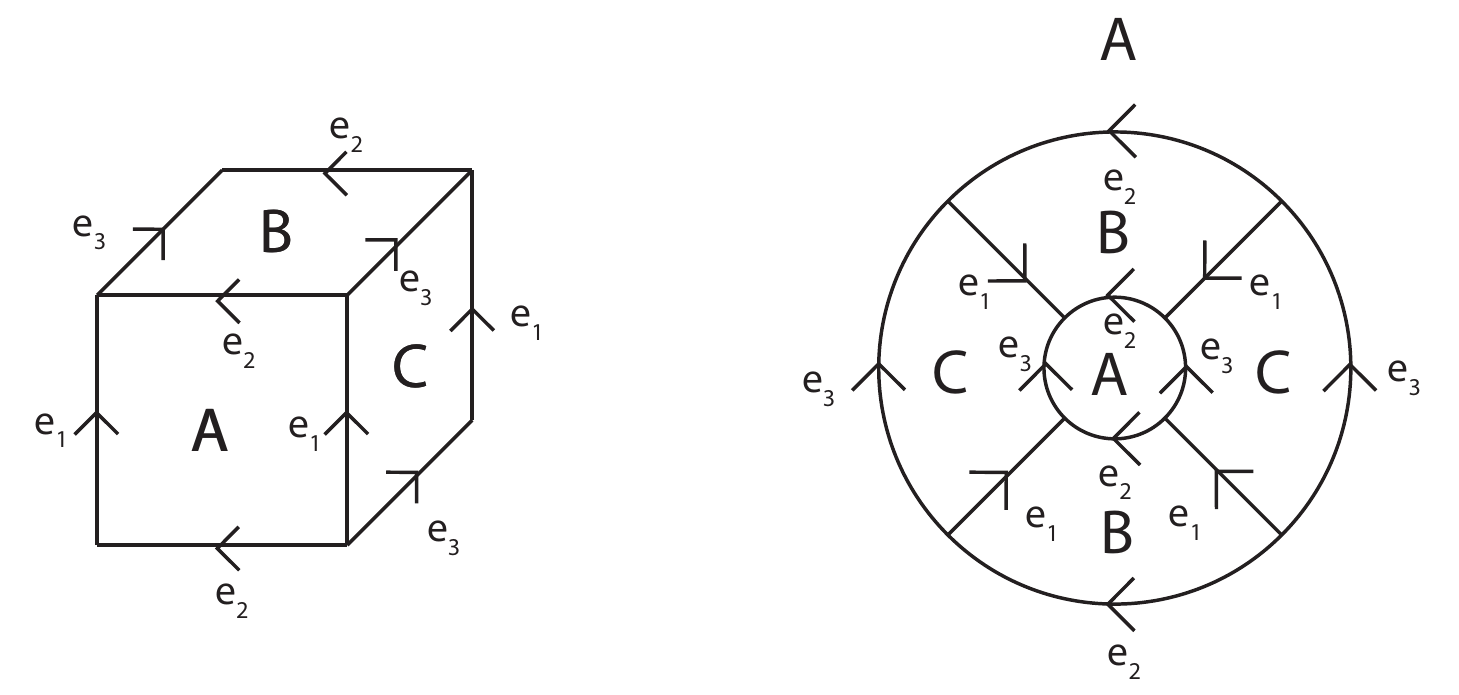}} \caption{The gluing map for $S^1\times S^1 \times S^1$.}
\label{3TorusGlue}
\end{center}
\end{figure}

The first geometry we examine is Euclidean 3-space. This is the most familiar geometry, and it is the nicest geometry for visualization. We use $S^1\times S^1 \times S^1$ as a representative manifold. A fundamental domain for $S^1\times S^1 \times S^1$ is the cube, with opposite faces identified, as shown in Figure \ref{3TorusGlue}. We will look at this example in great detail. As a preview, the subdivision rule we will obtain is illustrated in Figure \ref{CircleTorus}.

The universal cover of $S^1\times S^1 \times S^1$ is $\R^3$, and copies of this cubical fundamental domain tile three-space.  In this tiling, four cubes come together at an edge.  We describe this by saying each edge has edge cycle length 4. In general, when edges of a fundamental domain are glued together by a map, the number of edges in the equivalence class of an edge $e$ is the \textbf{edge cycle length} of $e$.

\begin{figure}
\begin{center}
\scalebox{.5}{\includegraphics{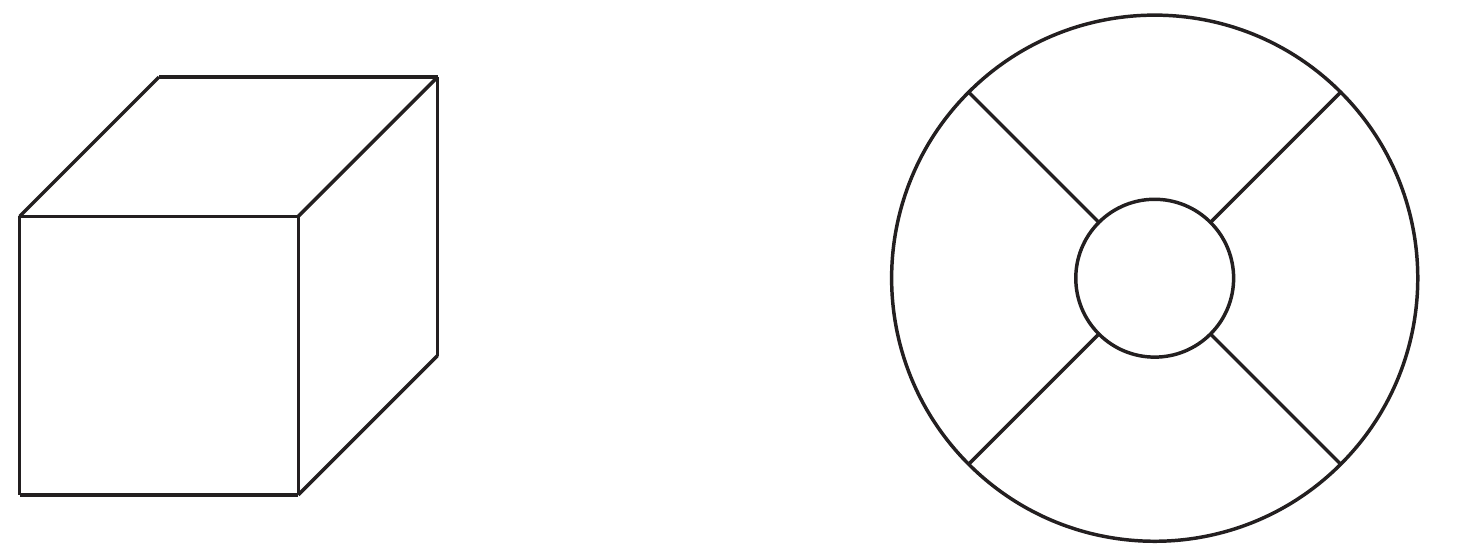}}
\caption{S(1)} \label{SOne}
\end{center}
\end{figure}

To find a subdivision rule, let's begin by constructing $S(n)$ for $S^1\times S^1 \times S^1$.  Recall that $B(n)$ is the `ball' of fundamental domains of distance at most $n$ from the identity in the word metric, and that $S(n)$ is the boundary of $B(n)$ ($B(n)$ may not always be a topological ball; closed $S^2\times\mathbb{R}$ manifolds are counterexamples). $S(1)$ is the projection of the cube. See Figure \ref{SOne}.

\begin{figure}
\begin{center}
\scalebox{.8}{\includegraphics{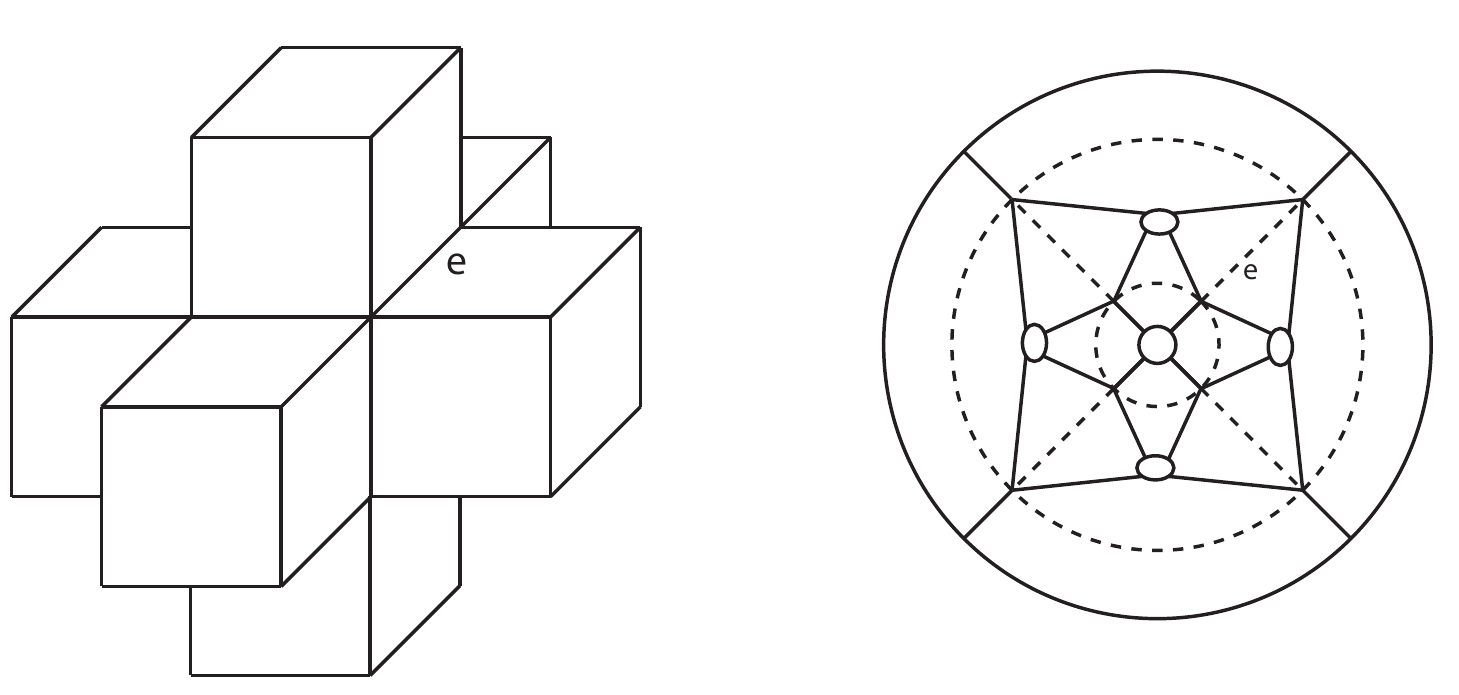}}
\caption{S(2)}\label{STwo}
\end{center}
\end{figure}

Adding a cube to each face, we get $S(2)$. See Figure \ref{STwo}. Dotted lines represent corners, or edges common to three copies of the fundamental domain.

This gluing is symmetric, so every face in $S(1)$ is replaced in the same way when creating $S(2)$. That is, each square face in $S(1)$ has a new cube glued on; the five remaining unglued faces of the cube form a disk with a new cell structure replacing the old square disk.

\begin{figure}
\begin{center}
\scalebox{.65}{\includegraphics{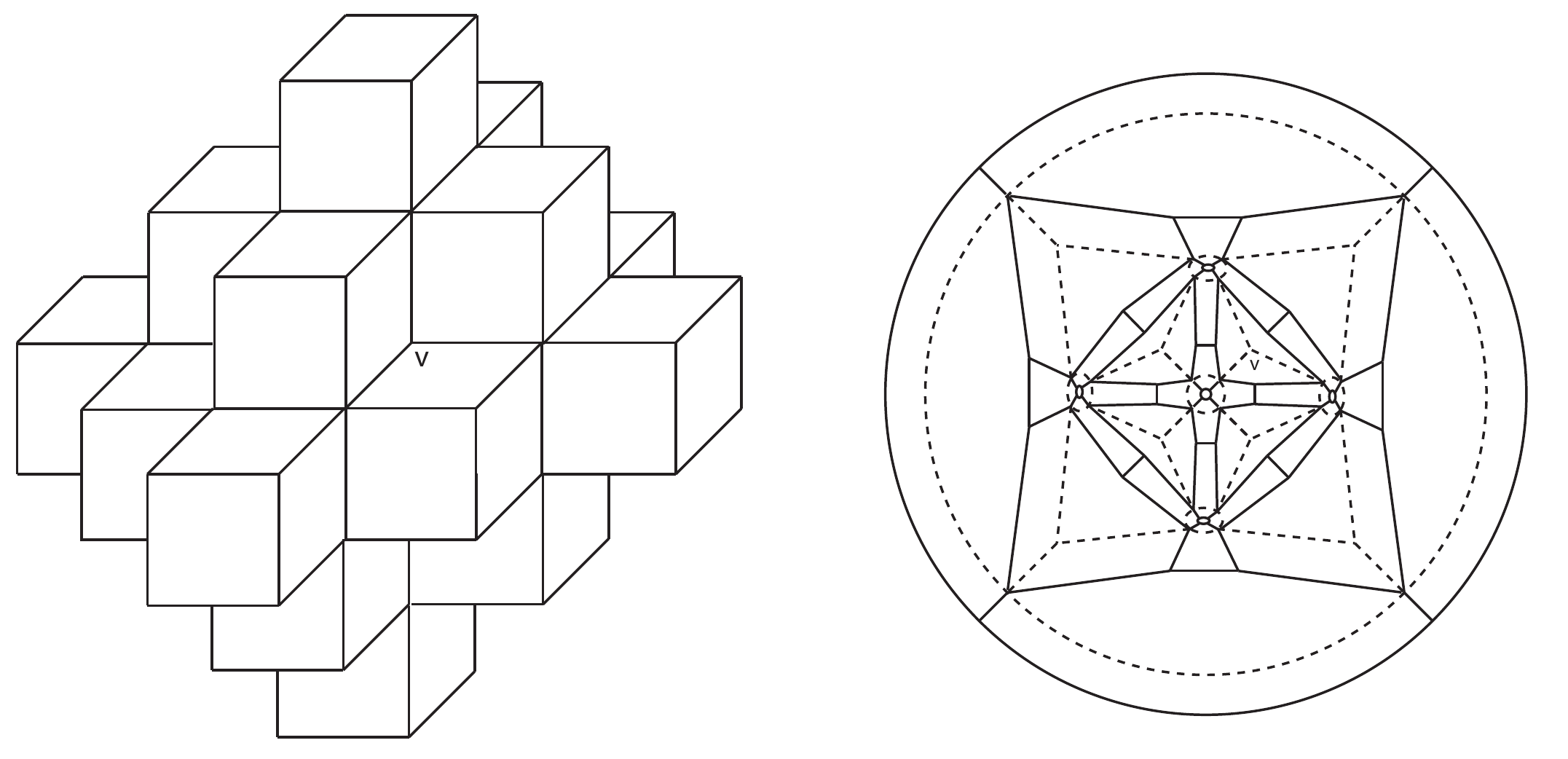}}
\caption{S(3)}\label{SThree}
\end{center}
\end{figure}

\begin{figure}
\begin{center}
\scalebox{.4}{\includegraphics{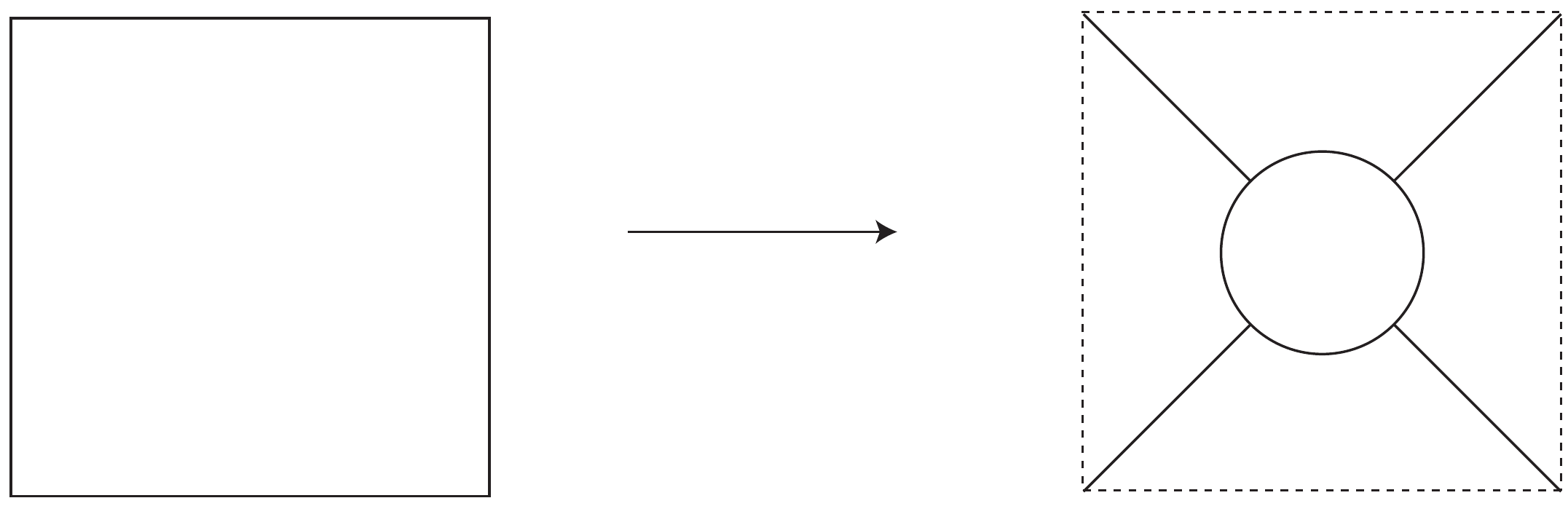}}
\caption{The replacement rule for a single face.}\label{TorusSquareSub}
\end{center}
\end{figure}

$S(3)$ is more complicated (see Figure \ref{SThree}). Notice that, in creating $S(3)$, every face in $S(2)$ without dotted edges is replaced just as the
original faces of $S(1)$ were. See Figure \ref{TorusSquareSub}.

\begin{figure}
\begin{center}
\scalebox{.4}{\includegraphics{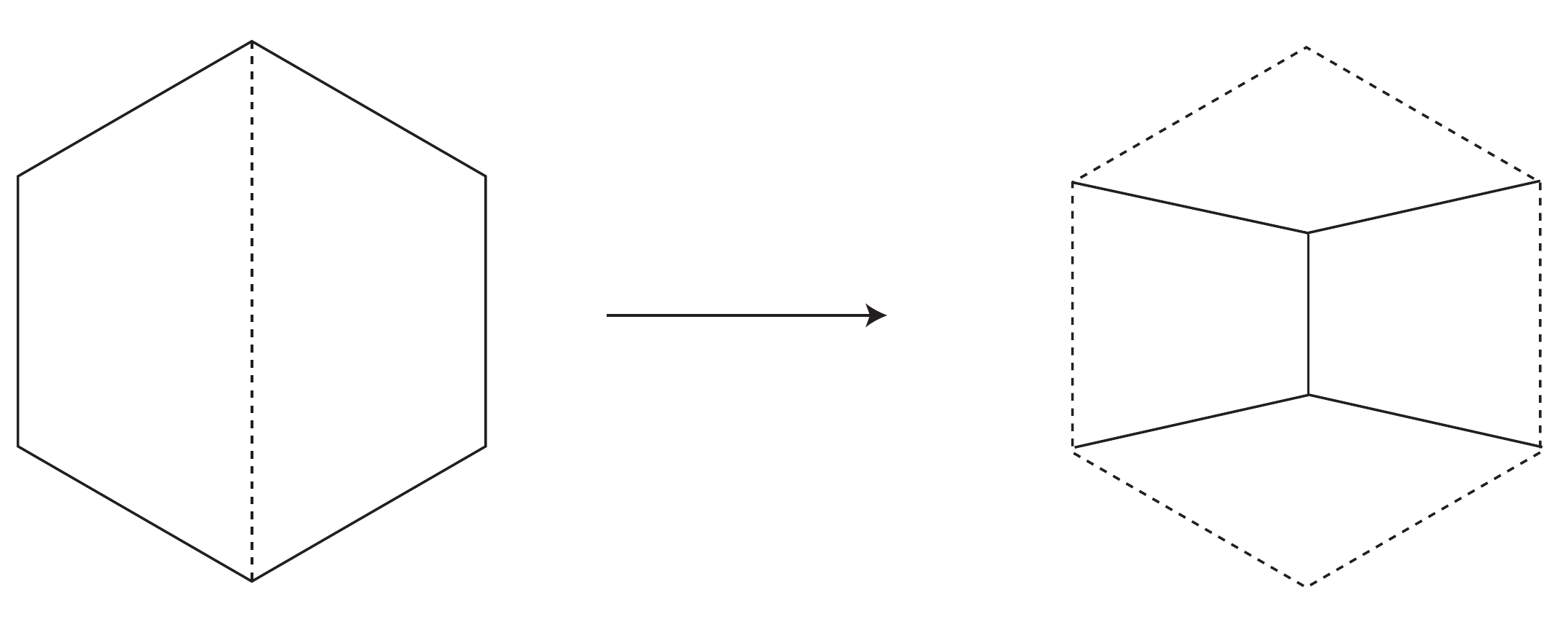}}
\caption{The replacement rule for a pair of faces.}\label{TorusPairSub}
\end{center}
\end{figure}

However, pairs of faces which share a dotted edge are replaced by a single polyhedron as in Figure \ref{TorusPairSub}. We will call the dotted edges, which correspond to corners, \textbf{loaded edges}. In general, a loaded edge in $S(n)$ is an edge that intersects (edge cycle length of $e)-1$ copies of the fundamental domain in $B(n)$.

\begin{figure}
\begin{center}
\scalebox{.6}{\includegraphics{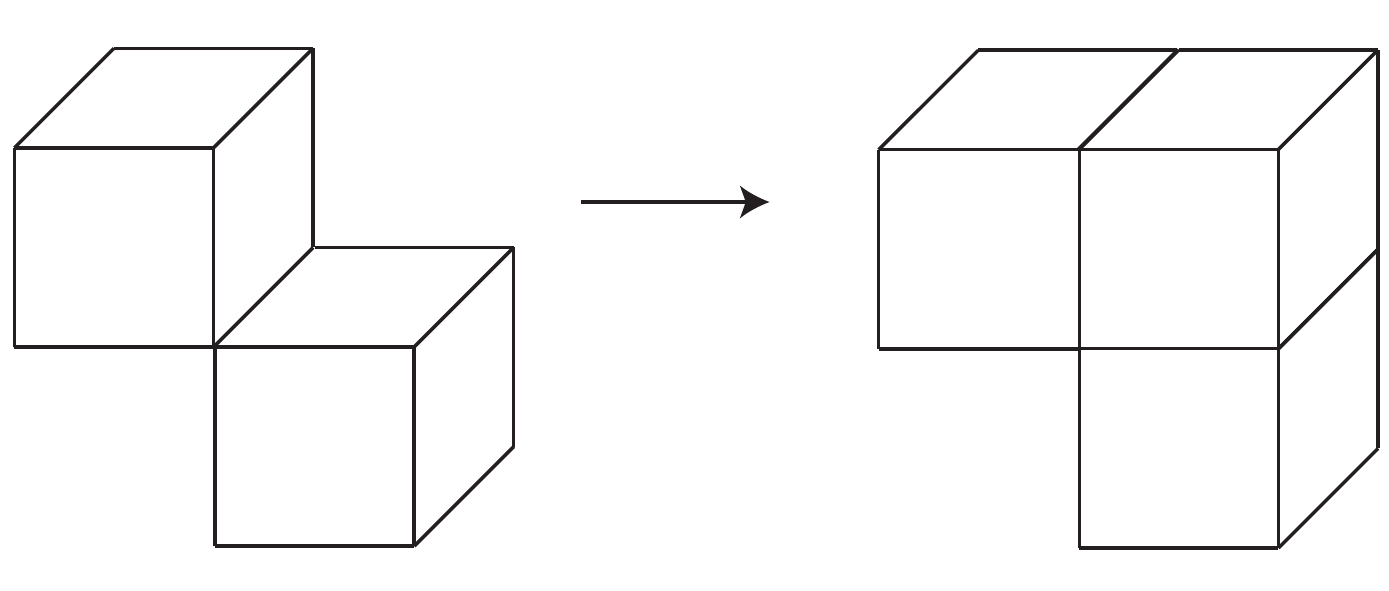}}
\caption{The 3D version of Figure \ref{TorusPairSub}.}\label{TorusCorner}
\end{center}
\end{figure}

Notice the loaded edge disappears in $S(3)$.  This corresponds to a corner in $S(2)$ being covered up by a single polyhedron (see Figure \ref{TorusCorner}). We glue a single polyhedron onto two faces because, in the universal cover, every edge should touch four cubes or copies of the fundamental domain. When an edge is loaded, it already intersects three fundamental domains, and so only one more can intersect that edge. Thus, a single cube must glue onto both faces.

\begin{figure}
\begin{center}
\scalebox{.4}{\includegraphics{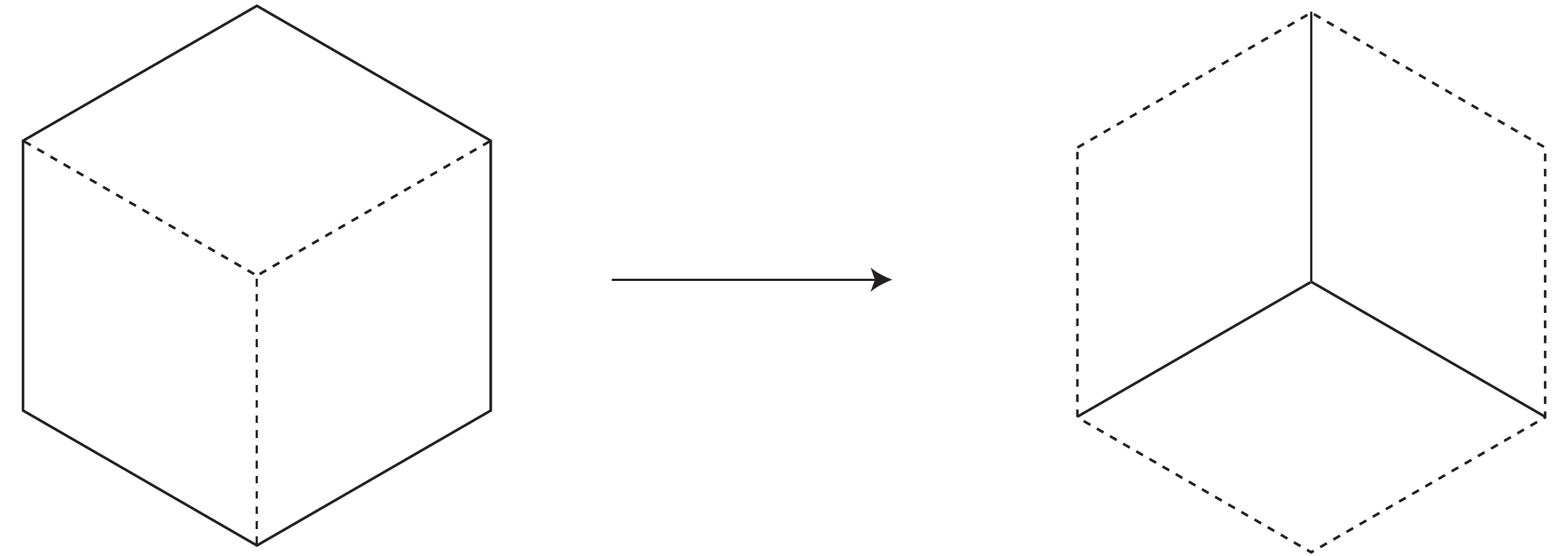}}
\caption{The replacement rule for three faces.}\label{TorusTripleSub}
\end{center}
\end{figure}

A third situation occurs in going from $S(3)$ to $S(4)$, where three loaded lines converge at a single vertex.  We will call this a \textbf{loaded vertex}.  In this case, a single polyhedron covers up all three, as shown in Figure \ref{TorusTripleSub}. In general, a vertex in $S(n)$ is a loaded vertex when all edges coming into it are loaded. A loaded vertex in $S(n)$ is always covered in $S(n+1)$.

\begin{figure}
\begin{center}
\scalebox{.9}{\includegraphics{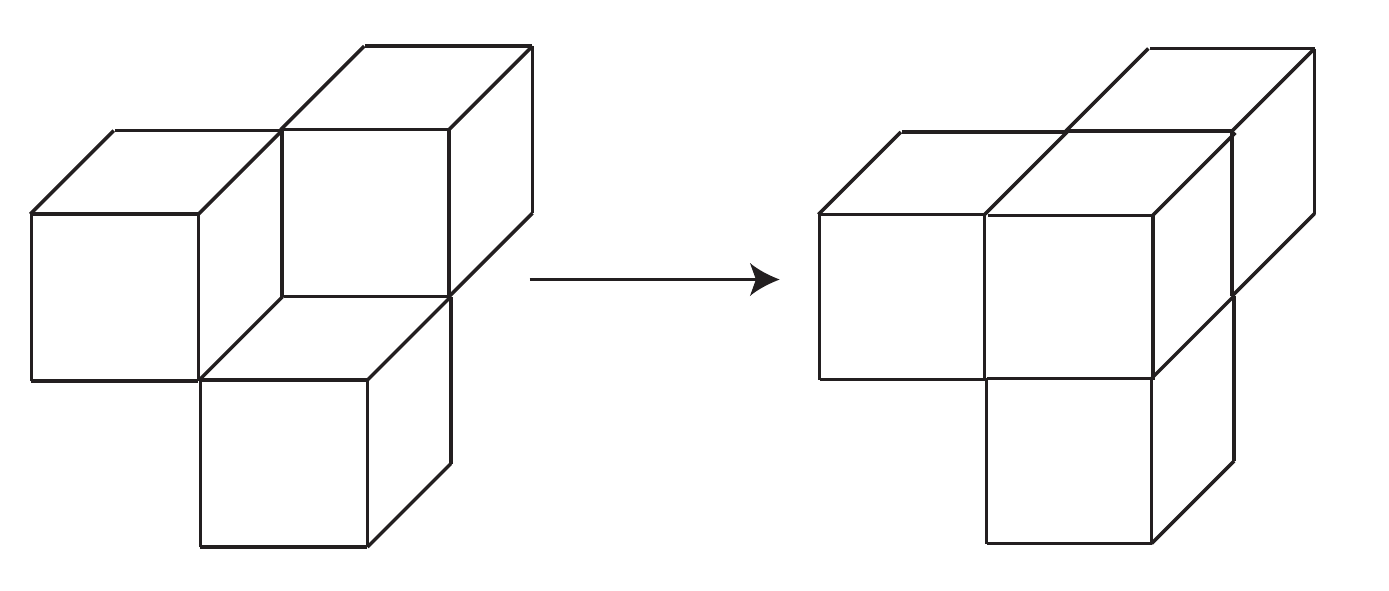}}
\caption{The 3D version of Figure \ref{TorusTripleSub}}\label{TorusBigCorner}
\end{center}
\end{figure}

This corresponds to the situation in Figure \ref{TorusBigCorner}.

Notice now that in $S(3)$, there are no new combinations of tiles; unloaded faces, two faces sharing a loaded edge and three faces sharing a loaded vertex are all that happen. It is clear that these situations are all that will ever happen. Notice that this is not a subdivision rule; edges are created, disappear, reappear, etc. However, we do have a replacement rule, as we know how to replace every local combination of faces that appear at every stage. This can be turned into a subdivision rule by eliminating the disappearance of cells. To do this, we add new edges at every stage.
 
\begin{figure}
\begin{center}
\scalebox{.6}{\includegraphics{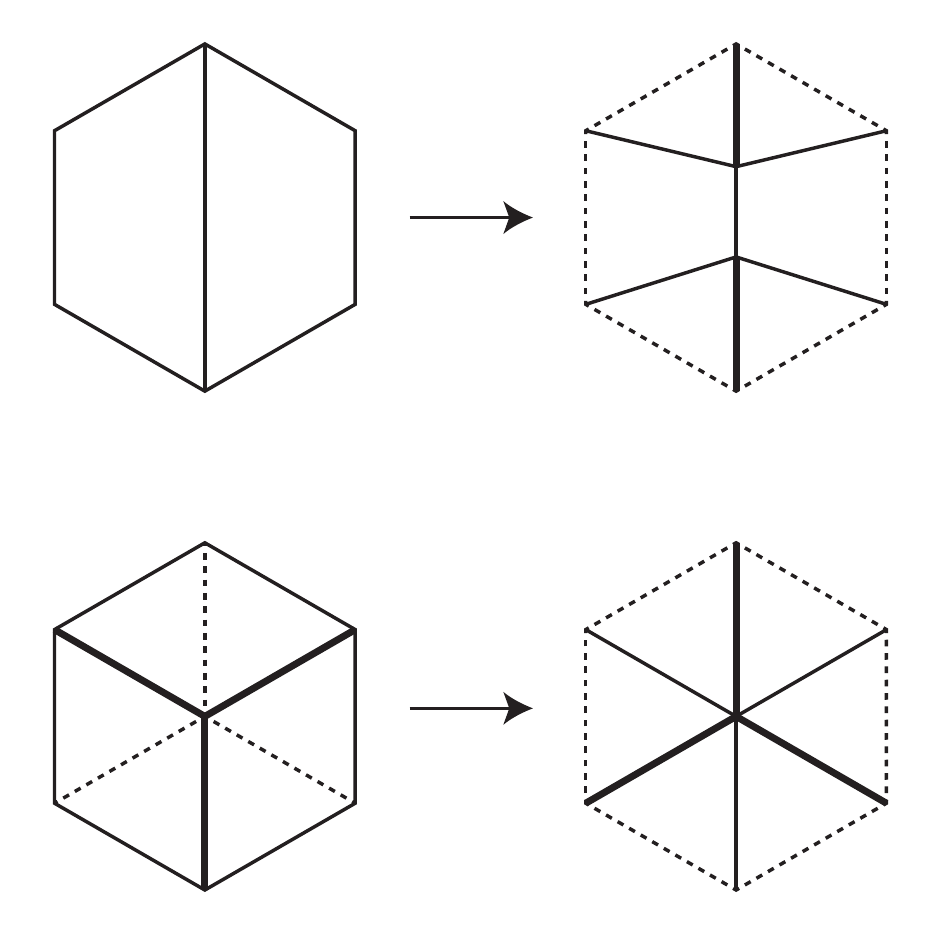}}
\caption{We can add lines to the two `loaded' tile types to get a self-consistent subdivision rule.}\label{TorusAddedLines}
\end{center}
\end{figure}

For instance, in Figure \ref{TorusPairSub}, the center line between two squares disappears when we glue on the new cube. However, if we add a line to the new cell structure (as shown in the top half of Figure \ref{TorusAddedLines}), then the new cell structure contains the old cell structure as a subset. Thus, we have a subdivision rule.

However, this divides the top and bottom squares (which are in $S(n+1)$, if the loaded pair they replaced was in $S(n)$) into two triangles each; each of these surround a loaded vertex, so we have to change the replacement rule for three squares surrounding a loaded vertex, (call the three faces a \textbf{loaded star}), since each square now has an extra line in $S(n+1)$. But notice that adding the three lines in to the loaded star on the left of Figure \ref{TorusTripleSub} gives us a hexagon divided into six `pie slices'. If we add similar lines bisecting the loaded star (as seen in Figure \ref{TorusAddedLines}), we again get a hexagon divided into six triangles; thus, the subdivision on each triangle in that hexagon is just the identity.  Adding these later 3 lines causes no new situations in $S(n+2)$, because each square with an added line is part of a loaded star. Thus, all loaded stars in all stages will have 3 extra lines.

\begin{figure}
\begin{center}
\scalebox{.6}{\includegraphics{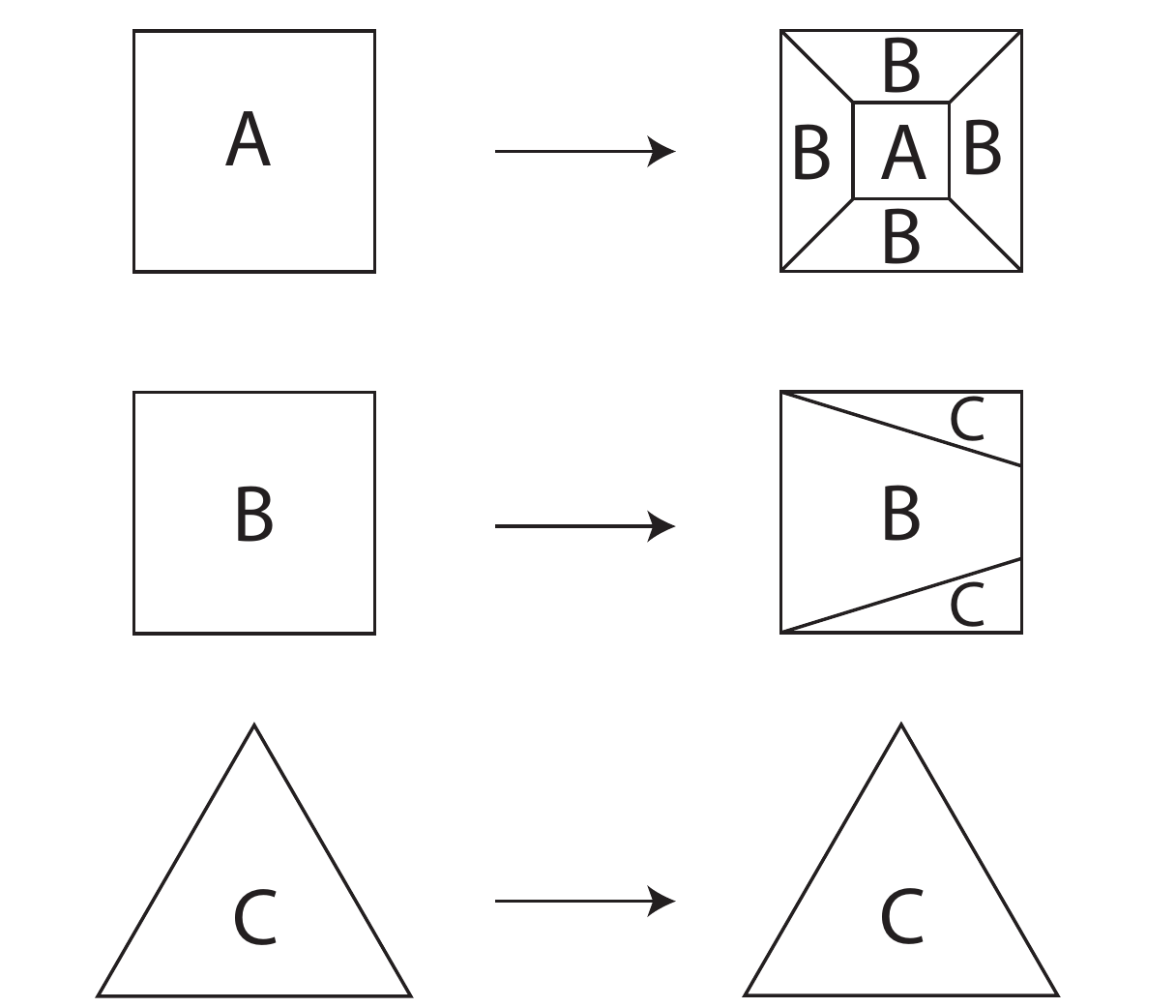}}
\caption{The replacement rule for the 3-dimensional torus.}\label{TorusSubs}
\end{center}
\end{figure}

We summarize this in Figure \ref{TorusSubs}. The circle packed pictures are shown on page \pageref{CircleTorus}. These circle packed pictures only display the subdivision rule combinatorially; the circle packed pictures are not isometrically subsets of each other, because this subdivision rule is not conformal. The connection between circle packings and conformality is explained in \cite{French}. Notice the similarity of this subdivision rule to the subdivision rule for the Hopf link (a Euclidean knot) from \cite{linksubs} (see Figure \ref{HopfBound}).

\section{\texorpdfstring{$\mathbb{H}^2\times \R$ geometry: $N\times S^1$}{H\texttwosuperior\ x\ R\ geometry:\ N\ x\ S\textonesuperior}}
\label{H2RSection}

In this section, we study the product geometry $\mathbb{H}^2 \times \mathbb{R}$, with example manifold $N \times \mathbb{S}^1$.

\begin{figure}
\begin{center}
\scalebox{.8}{\includegraphics{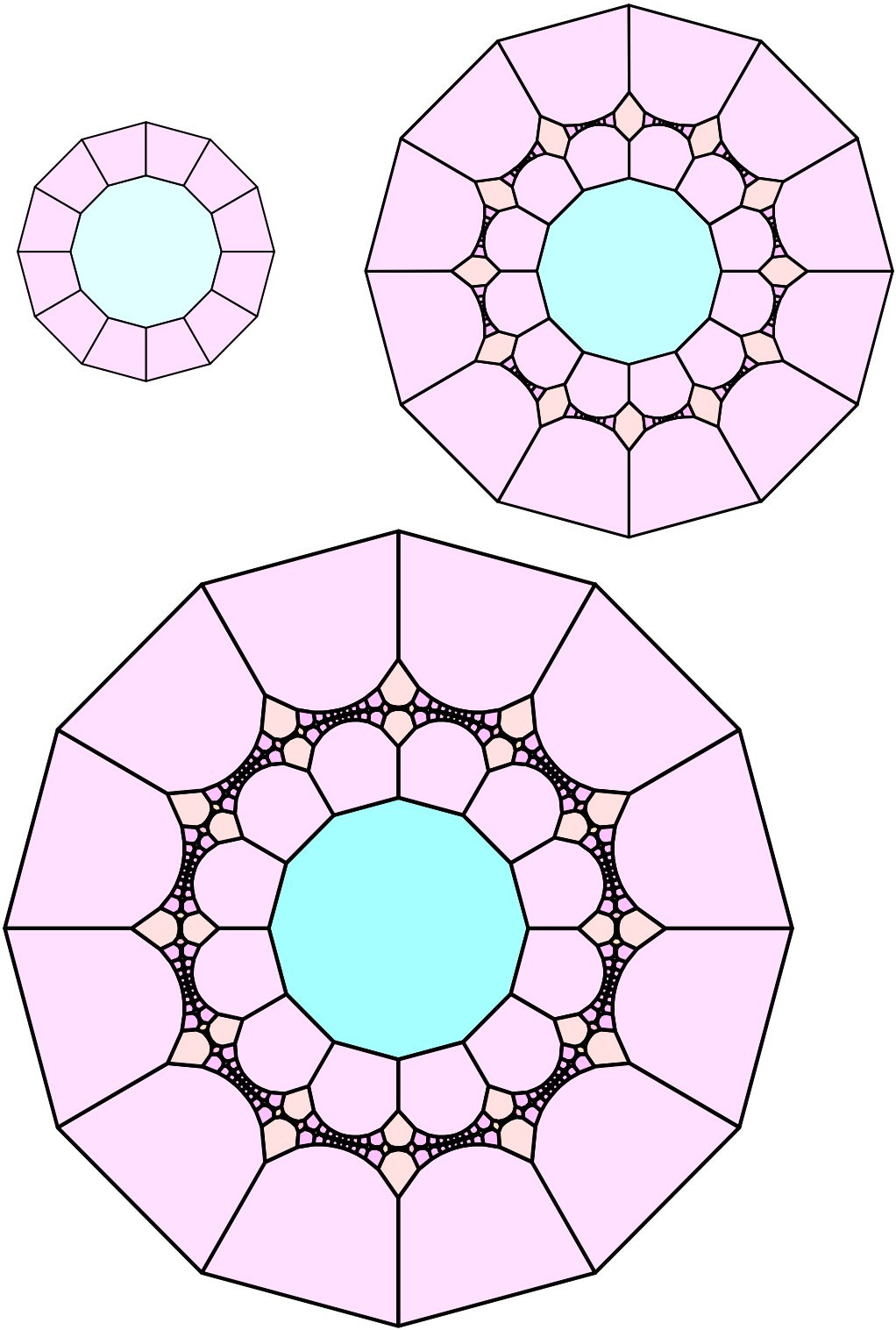}} \caption{$R(X)$, $R^2(X)$, and $R^3(X)$, where $X$ is the top and sides of a fundamental domain for $N\times S^1$.}
\label{CircleProduct}
\end{center}
\end{figure}

Here we are using $N$ to represent the non-orientable surface of Euler characteristic $-1$, i.e. the connected sum of three projective planes.  We chose this surface as our example because every other hyperbolic surface can be pieced together from it.  In this sense, it is the smallest hyperbolic surface.  We also chose this example because it introduces the notion of `fragile edges', which come from manifolds with edges of odd cycle length. The circle packed picture of the subdivision rule we will get is shown in Figure \ref{CircleProduct}.

\begin{figure}
\begin{center}
\scalebox{1.2}{\includegraphics{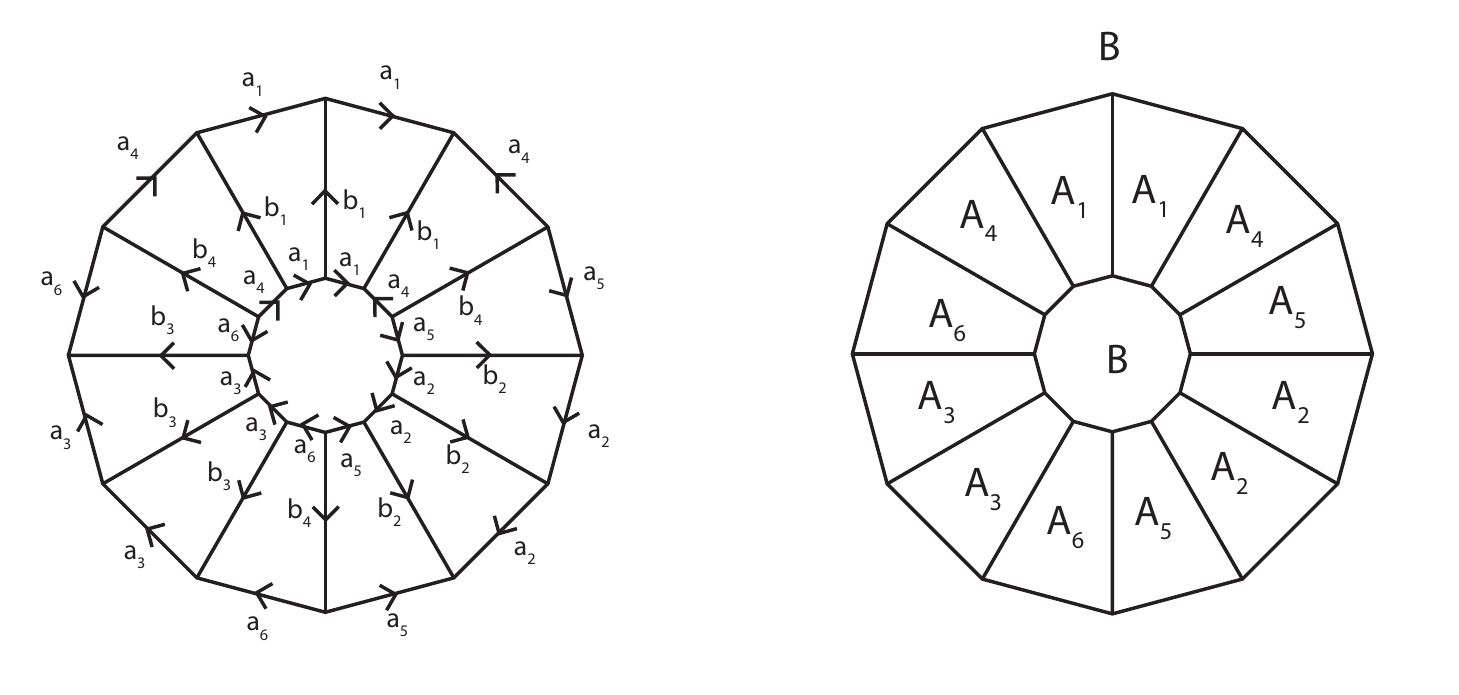}} \caption{The gluings for our $N\times S^1$ manifold.}
\label{DodoGlue}
\end{center}
\end{figure}

The 3-manifold $N\times S^1$ has a fundamental domain that is a dodecagonal prism. Figure \ref{DodoGlue} shows the necessary gluings.

Notice that each $a_i$ edge has edge cycle length 3, and each $b_i$ edge has edge cycle length 4.  This particular gluing was chosen to give these cycle lengths.  Some of the face-gluing maps are orientation reversing, and others are orientation preserving.  But, in creating $S(n)$, changing the orientation won't change the combinatorial structure of the fundamental domains.  This is because the prism has symmetry group $D_{12}\times Z_2$ and preserves its shape under reflection.

We again let $B(1)$ be a single fundamental domain with $S(1)$ its boundary.  Now that we are dealing with hyperbolic space, the balls $B(n)$ and spheres $S(n)$ are more difficult to imagine.  But we can describe combinatorially how each face is replaced, as we did for $S^1 \times S^1 \times S^1$ earlier.

\begin{figure}
\begin{center}
\scalebox{.7}{\includegraphics{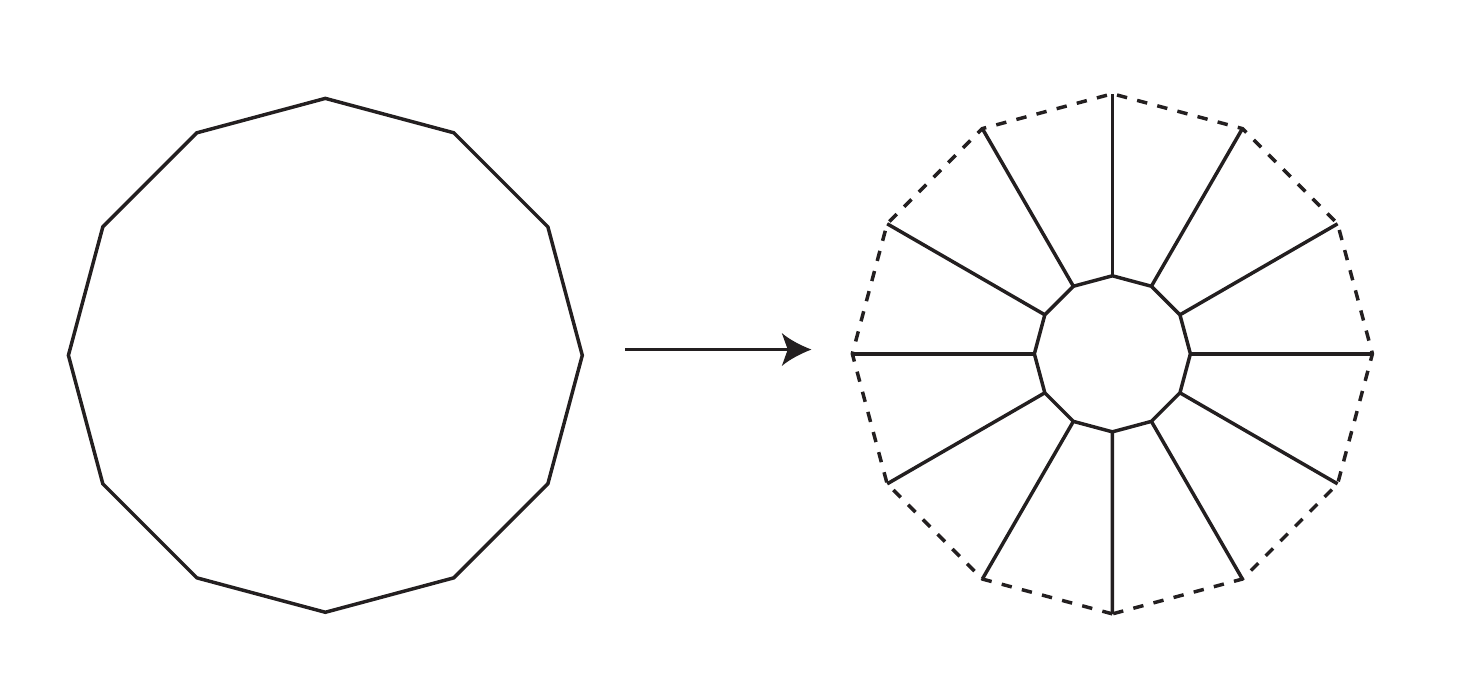}}
\caption{Type B face.}\label{DodoBSub}
\end{center}
\end{figure}

\begin{figure}
\begin{center}
\scalebox{.7}{\includegraphics{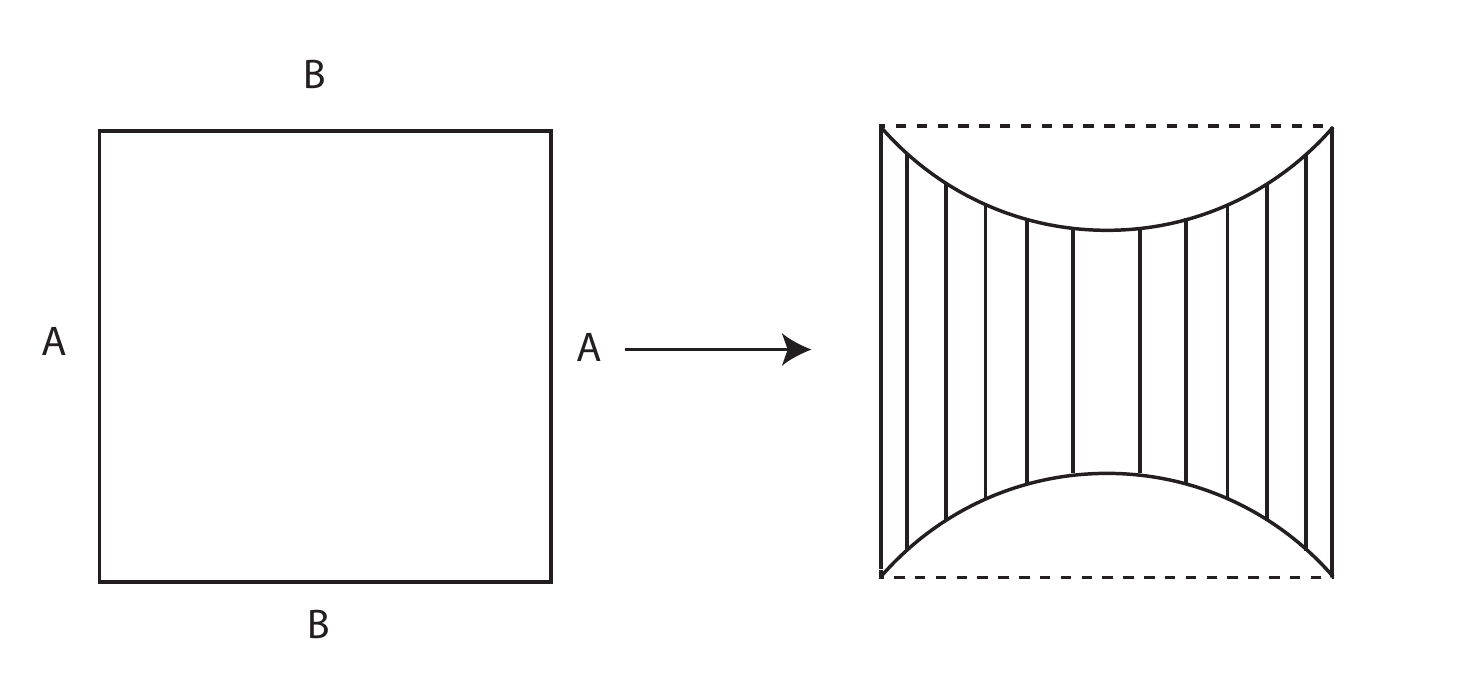}}
\caption{Type A face.}\label{DodoASub}
\end{center}
\end{figure}

We call all dodecagonal faces on $S(1)$ type B, and all square faces type A. Faces of type B are replaced as in Figure \ref{DodoBSub}. Similarly, type A faces are replaced as in Figure \ref{DodoASub}.

\begin{figure}
\begin{center}
\scalebox{.8}{\includegraphics{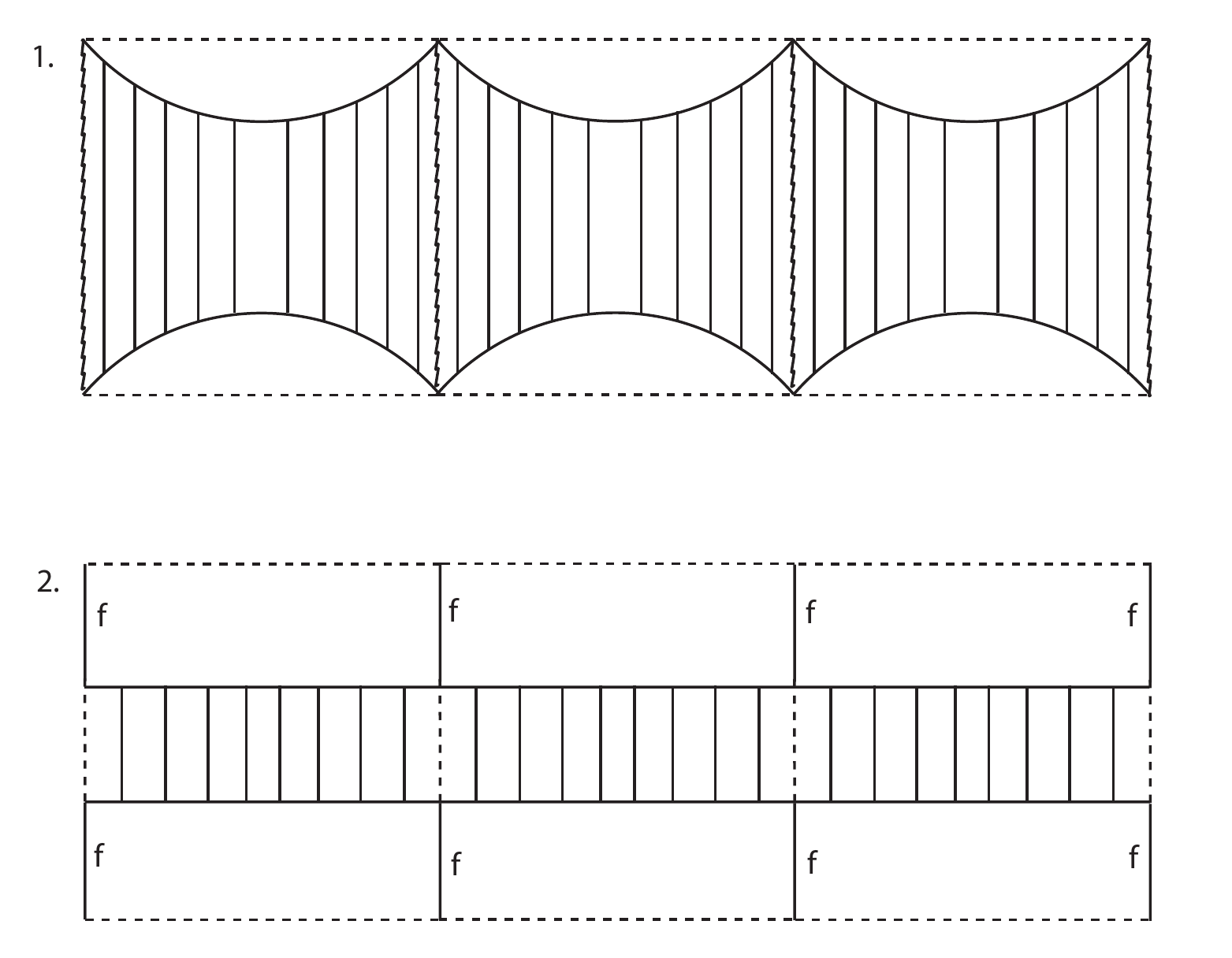}} \caption[Collapsing of edges.]{1. The squiggly lines represent edges about to collapse.  2. After collapse.  The edges marked with f are now `fragile'. All square faces without dotted edges or f's act like type A faces.}
\label{DodoCollapse}
\end{center}
\end{figure}

However, here we have a problem.  The $a_i$ edges between A faces have edge cycle length 3.  Thus, only three polyhedra can intersect that edge in the universal cover.  In $B(2)$, three polyhedra already intersect $a_i$: the original polyhedron, and the two polyhedra that were glued onto the faces on either side of the edges in $S(2)$.  There are now two unglued faces on either side of any given $a_i$ edge; if we glue new polyhedra to both faces, we'll have too many polyhedra intersecting that edge.  Thus, since these two faces must be glued to something,  and there is no other choice that will keep our $B(2)$ simply connected, we must glue them to each other, as shown in Figure \ref{DodoCollapse}.

The essential point here was that only two more polyhedra could intersect the edge, causing the faces bordering that edge to collapse.  Any edge $e$ with edge cycle length $L$ that already intersects $L-2$ polyhedra will be called \textbf{fragile}.  Edges can be become fragile over time, just as they can become loaded. Note that burdened edges are different from fragile edges.

In our replacement of the A tiles above, all of the visible edges of the faces that collapsed together now intersect two polyhedra.  Thus, edge cycle length 3 edges are loaded, and edge cycle length 4 edges are fragile.

These replacements have created new types of faces to deal with.  Specifically, we now have: loaded pairs of B's and A's, where the B portion has two fragile edges, as well as loaded pairs of A's.

\begin{figure}
\begin{center}
\scalebox{.7}{\includegraphics{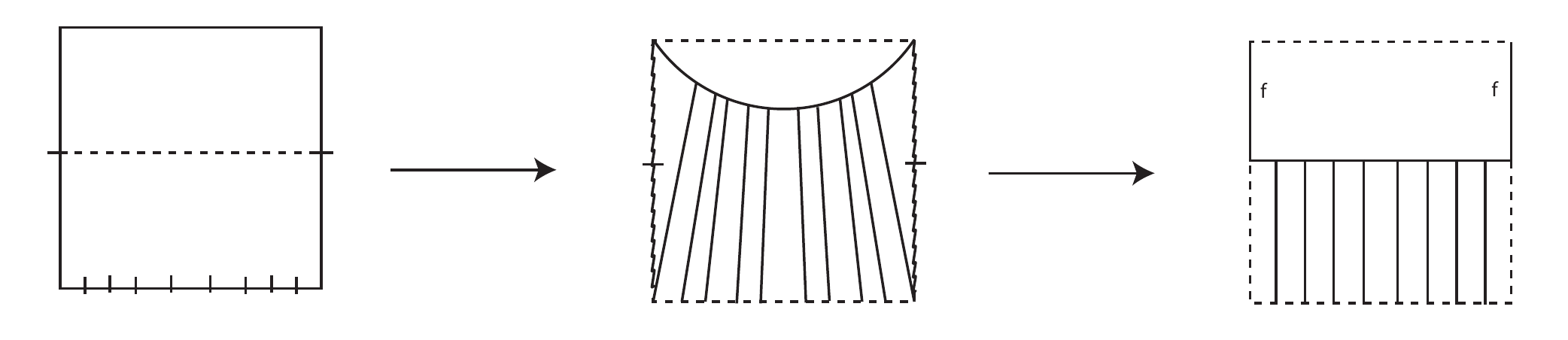}} \caption{The replacement rule for an A/B pair. All squares on the bottome besides the two corners are parts of A/B pairs. The top is the B part of an A/B pair. The corner tiles are parts of A/A/B triples.}
\label{DodoABPairRule}
\end{center}
\end{figure}

The A/B pairs are replaced as in Figure \ref{DodoABPairRule}.

Notice that the leftmost and rightmost regions collapse along two fragile edges each.  It is possible that this would cause confusion, as each of those faces must be identified to both faces that it borders over those two edges.  However, looking carefully, we see that A/B pairs border only other A/B pairs on either side, so that the collapsing is well-defined.
Note that this covers up an old vertex and adds a new one.  For convenience, we've placed the new vertex directly over the old.

\begin{figure}
\begin{center}
\scalebox{.75}{\includegraphics{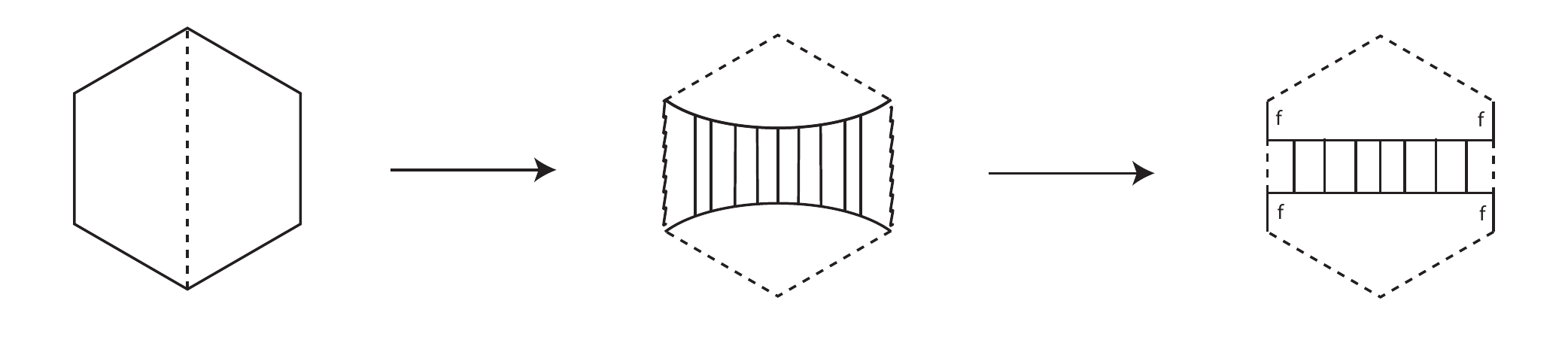}} \caption{The replacement rule for an A/A pair. The top and bottome tiles are B parts of A/A/B triples.}
\label{DodoAAPairRule}
\end{center}
\end{figure}

Pairs of A tiles behave as in Figure \ref{DodoAAPairRule}.

\begin{figure}
\begin{center}
\scalebox{.7}{\includegraphics{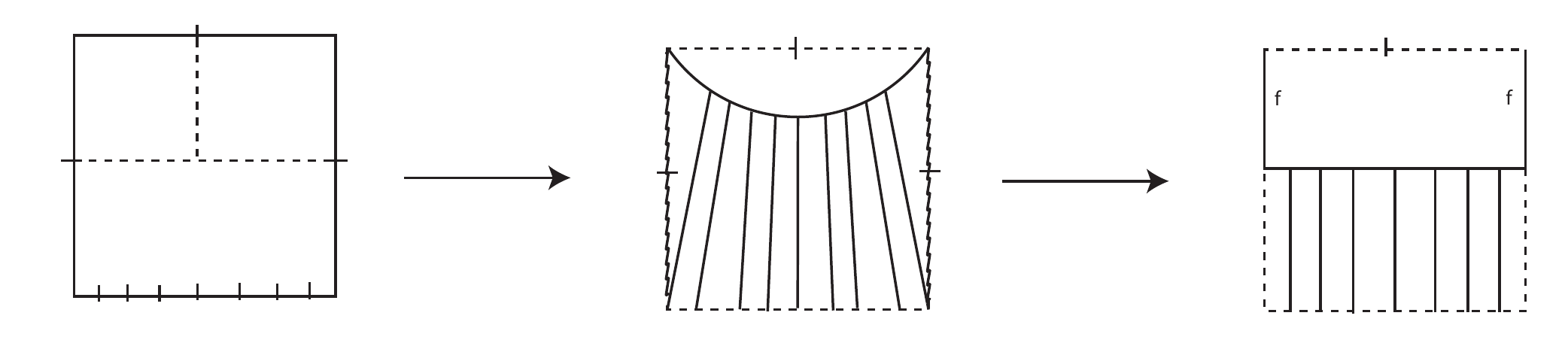}} \caption{The replacement rule for an A/A/B triple. The dodecagon on top is the B part of an A/A/B triple. The bottom squares are all A tiles in A/B pairs, except the corners, which are A tiles of A/A/B triples.}
\label{DodoTripleRule}
\end{center}
\end{figure}

Together, both kinds of pairs create loaded triples consisting of two A's and a B, where the B has two fragile edges.  These triples are subdivided as in Figure \ref{DodoTripleRule}.

Again, the leftmost and rightmost regions collapse over two edges each.  No new types of faces have occurred, so we have found all of the replacement rules.

\begin{figure}
\begin{center}
\scalebox{.8}{\includegraphics{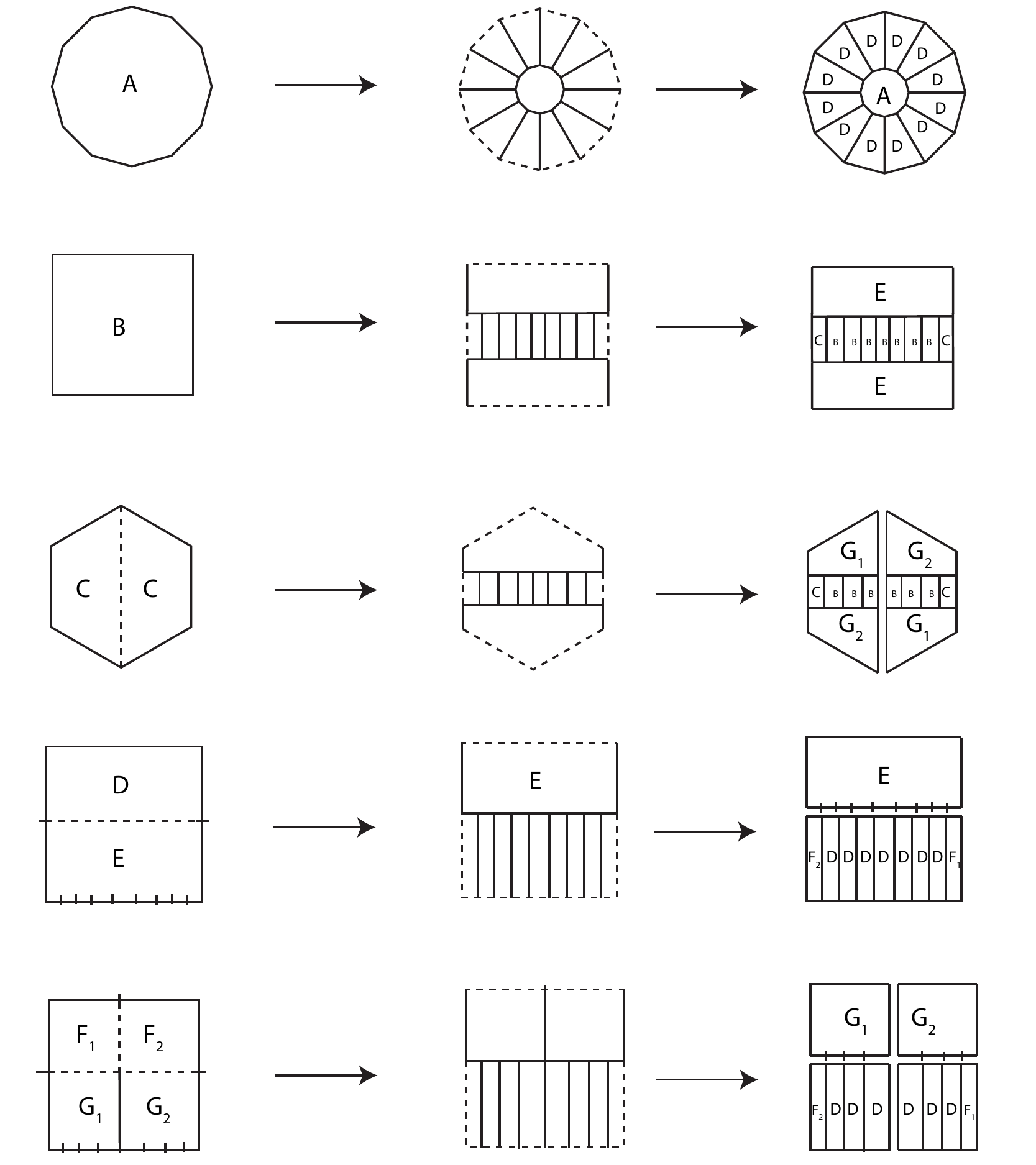}} \caption{The subdivision rule for M.}
\label{DodoSubs}
\end{center}
\end{figure}

To convert this to a subdivision rule, we must add lines to make the cell structure of $S(n+1)$ a refinement of the cell structure of $S(n)$. These lines are added as shown in Figure \ref{DodoSubs}.  We add one vertical line to A/A pairs after replacement to create C tiles. The cell structure of A/B tiles before refinement already imbeds into the replaced cell structure, and we add no lines. The lines we added to all A/A pairs means that A/A/B triples already have extra lines before replacement. We add a vertical line to the B portion of all A/A/B triples after replacement, and this is all that is needed. Because we added a line to all B tiles in A/A/B triples, and we added a line to all B tiles in A/A pairs, all B tiles that are part of A/A/B triples will consistently have a line added. All of this gives us a subdivision rule. The circle packings of this subdivision rule are shown in Figure \ref{CircleProduct}. This subdivision rule is very similar to the subdivision rule for the trefoil knot (and other 2-braid knots besides the Hopf link) as described in \cite{linksubs} and shown in Figure \ref{trefbound}.

\section{\texorpdfstring{$\mathbb{H}^3$ geometry: Hyperbolic dodecahedral space}{H\textthreesuperior\ geometry:\ Hyperbolic\ dodecahedral\ space}}
Subdivision rules were originally found for hyperbolic geometry. Cannon and Swenson's work \cite{hyperbolic} shows that closed manifolds with $\mathbb{H}^3$ geometry have subdivision rules in some sense. These subdivision rules are conformal; in particular, their combinatorial mesh goes to 0.

\begin{figure}
\begin{center}
\scalebox{.8}{\includegraphics{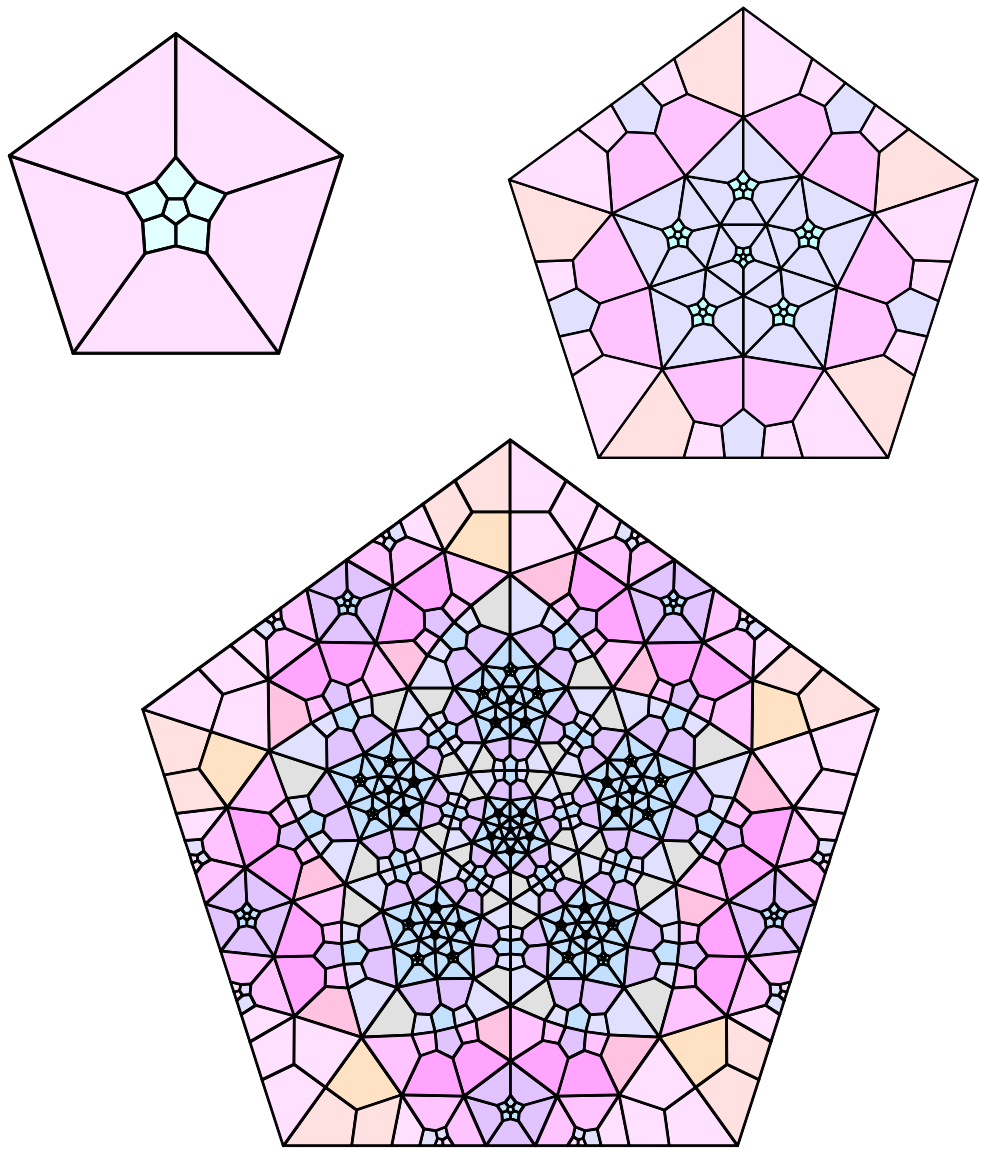}} \caption{$R(X)$, $R^2(X)$, and $R^3(X)$ for the subdivision rule of the right-angled dodechahedral orbifold, where $X$ is one face of the dodecahedral fundamental domain.}
\label{CircleHyperbolic}
\end{center}
\end{figure}

\begin{figure}
\begin{center}
\scalebox{.8}{\includegraphics{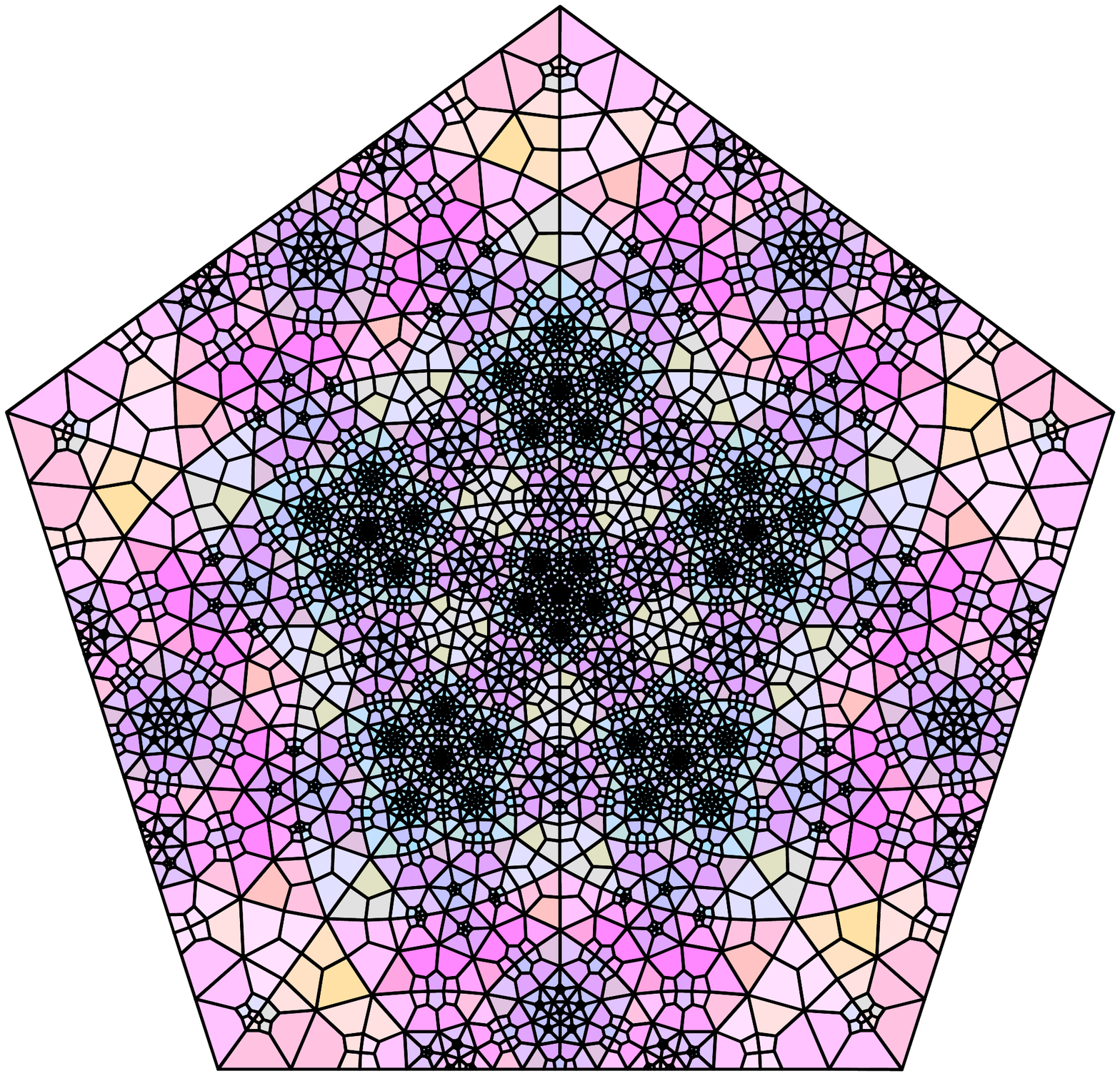}} \caption{$R^4(X)$ for the same manifold.}
\label{CircleHyperbolicBig}
\end{center}
\end{figure}

In \cite{myself2}, we found finite subdivision rules for all manifolds created by gluing right-angled hyperbolic polyhedra. One example is shown in Figures \ref{CircleHyperbolic} and \ref{CircleHyperbolicBig}. In general, the subdivision rules have mesh going to 0, and are the type of finite subdivision rules studied extensively by Cannon et. al., for instance in \cite{subdivision}. These subdivision rules are very similar in character to the subdivision rules of the hyperbolic alternating links found in \cite{linksubs}, which are shown in Figure \ref{Bringbound}.

\section{\texorpdfstring{$S^2\times \R$ geometry: $S^2\times S^1$}{S\texttwosuperior\ x\ R\ geometry:\ S\texttwosuperior\ x\ S\textonesuperior}}
Subdivision rules for $S^2\times \R$ manifolds are different than all others. The universal cover has two boundaries, which can be thought of as an `inner' and an `outer' sphere at infinity.

\begin{figure}
\begin{center}
\scalebox{.835}{\includegraphics{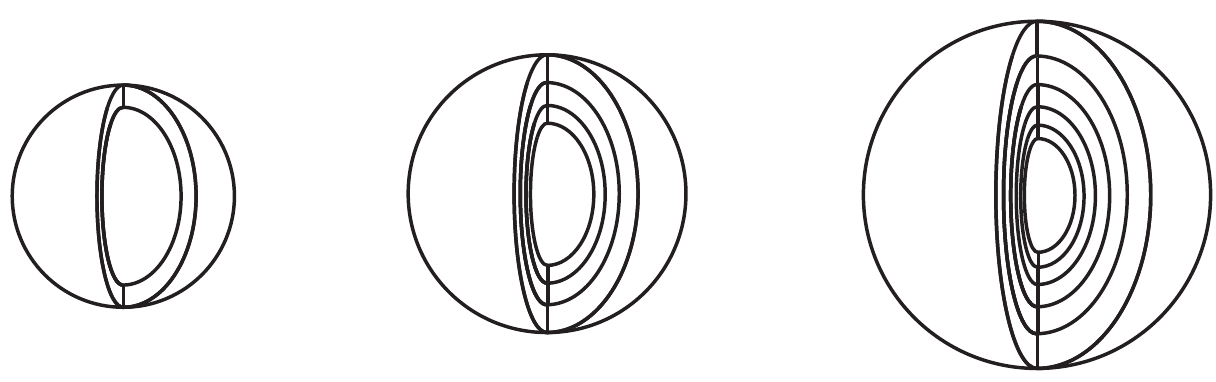}} \caption{The first
three stages in constructing the universal cover of $S^2\times
S^1$.} \label{S2Spheres}
\end{center}
\end{figure}

Consider $S^2\times S^1$.  This has as a fundamental domain
$S^2\times I$, a thickened sphere.  Constructing the universal cover
amounts to nesting these thickened spheres.  See Figure
\ref{S2Spheres}.

The figures for $S^2\widetilde{\times}S^1$ and $P^3\oplus P^3$ are
similar, as they also have a thickened sphere for their fundamental
domain.  $P^2\times S^1$ is only slightly more complicated.
In all of these cases, there are a fixed number of faces on each
boundary sphere at each sufficiently large stage of subdivision, giving us a sequence
of tilings $S(n)$ which is constant. The subdivision rule can be taken to be the identity on both spheres.

\section{\texorpdfstring{$\widetilde{SL_2}(\R)$ geometry: The unit tangent bundle of $N_{-1}$}{SL2(R)\ geometry:\ The\ unit\ tangent\ bundle\ of\ N}} \label{ProductMani}
Let's now consider the unit tangent bundle of $N$.  We choose this particular example of an $\tilde{SL_2}(\R)$ manifold for the same reason we chose $N\times S^1$ as an example for $H^2\times\R$ earlier: it is covered by all other unit tangent bundles of hyperbolic surfaces, and these are the best known manifolds corresponding to this geometry.  In fact, every $\widetilde{SL_2}(\R)$ manifold is finitely covered by the unit tangent bundle of a hyperbolic surface. For convenience, we will refer to the unit tangent bundle of $N$ as M throughout this section.

\begin{figure}
\begin{center}
\scalebox{1.2}{\includegraphics{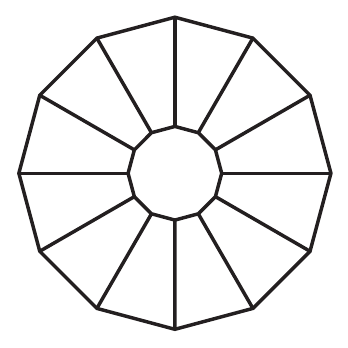}} \caption{The manifold
M after slicing.} \label{SlrStart}
\end{center}
\end{figure}

How can we find a polyhedral fundamental domain for M?  We can
project M onto $N$.  Slice $N$ as we did in creating
Figure \ref{DodoGlue} and slice M along the pre-images of
these edges.  What remains is the unit tangent bundle of a closed
disk.  This is necessarily the trivial bundle, and so we can slice
it along a fixed, `horizontal' surface to get the product of a
closed disk and the closed interval.  But this is just a 3-ball, and
so we have our polyhedron.  See Figure \ref{SlrStart}.

To get the manifold back, or to understand the combinatorics, we
need to understand how the edges and faces are glued.  The top is
glued to the bottom, but what about the square faces?

\begin{figure}
\begin{center}
\scalebox{.96}{\includegraphics{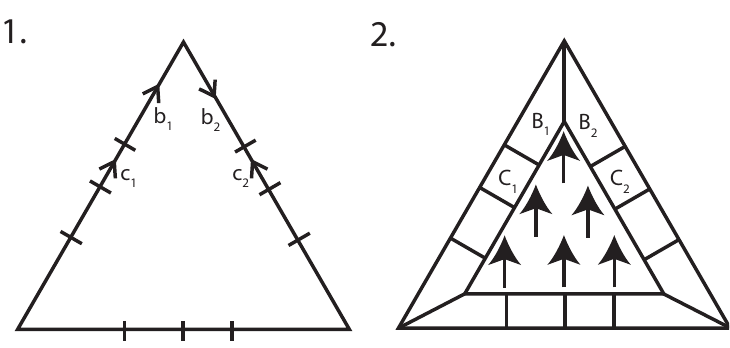}} \caption{The choice of
vector field for our horizontal slice.} \label{SlrField}
\end{center}
\end{figure}

Remember that we got a polyhedron by slicing the unit tangent bundle
of the disk along a horizontal surface.   A horizontal surface
corresponds to a fixed unit vector field $\alpha(x,y)$ on the disk.
Let's make the choice of vector field explicit, as shown in Figure
\ref{SlrField}. Notice that we have redrawn our dodecagon on as triangle.

Then we make a correspondence between the height of a point in the
polyhedron and the angle in the unit tangent bundle.  See Figure
\ref{SlrField}. If the top is 0 and the bottom is $2 \pi$ (viewing
the closed interval as $[0,2 \pi]$), then the angle at a point
$(x,y,z)$ is $\alpha(x,y)+z$.  This enables us to define the gluing
map.

\begin{figure}
\begin{center}
\scalebox{.8}{\includegraphics{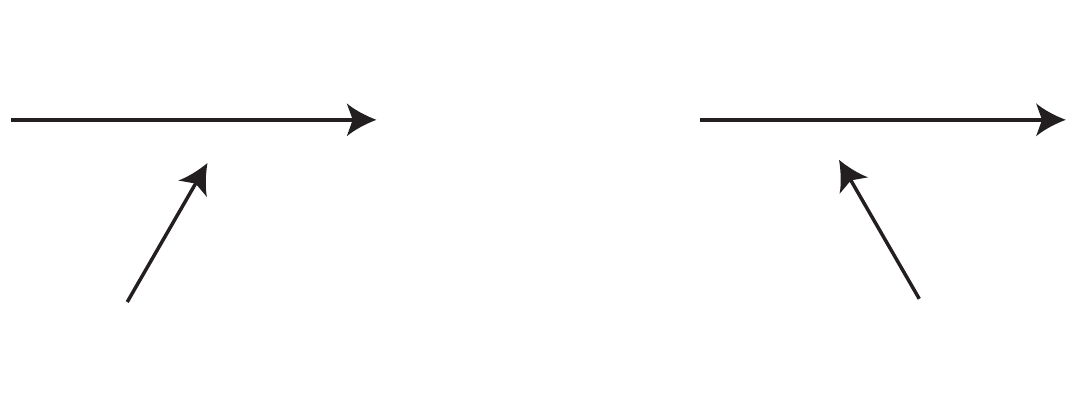}} \caption[The angles of the vector field.]{The angles
that the vector field makes with the edges $a_1$ and $a_2$ are
$\frac{\pi}{3}$ and $\frac{2 \pi}{3}$. } \label{SlrBAngles}
\end{center}
\end{figure}

\begin{figure}
\begin{center}
\scalebox{.8}{\includegraphics{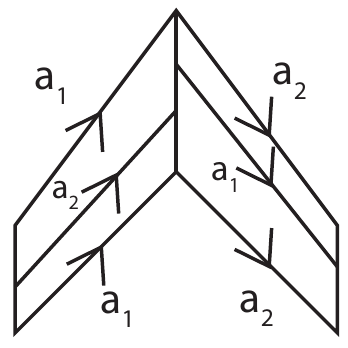}} \caption[The gluing diagram for an $\tilde{SL_2(R)}$ manifold]{The depth
in the z direction measures the change tangent direction from the
horizontal field at the top. Since the vector field at $a_1$ comes
in at an angle $\frac{\pi}{3}$ more than the vector field at $a_2$,
it gets mapped to a line at depth $\frac{\pi}{3}$.  Similarly, the
edge $a_2$ gets mapped to a depth of $\frac{5 \pi}{3}$.}
\label{SlrBGluing}
\end{center}
\end{figure}

Look at the faces $B_1$ and $B_2$.  They are identified together
because their projections in $N$ are edges that are identified.
Now look at the edge $a_1$.  Relative to the face $B_1$, the unit vector
field comes in at an angle of $\frac{\pi}{3}$, at the top of $B_1$
(i.e. at the edge $a_1$) , and at an angle of $\frac{2 \pi}{3}$ at
the top of $B_2$ (i.e. the edge $a_2$). See Figure \ref{SlrBAngles}.
This means that $a_1$ gets sent to a line on $B_2$ of depth
$\frac{\pi}{3}$. Similarly, $a_2$ gets sent to a line on $B_1$ with
depth $\frac{5 \pi}{3}$, so the faces $B_1$ and $B_2$ get mapped together as
shown in Figure \ref{SlrBGluing}.

\begin{figure}
\begin{center}
\scalebox{.8}{\includegraphics{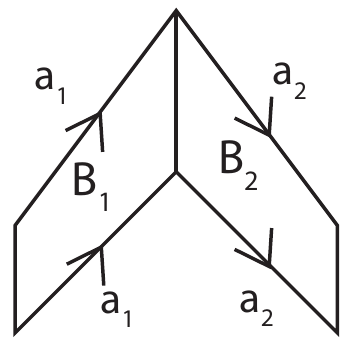}} \caption{The faces
$B_1$ and $B_2$ and the edges $a_1$ and $a_2$. } \label{SlrBZoomup}
\end{center}
\end{figure}

\begin{figure}
\begin{center}
\scalebox{.8}{\includegraphics{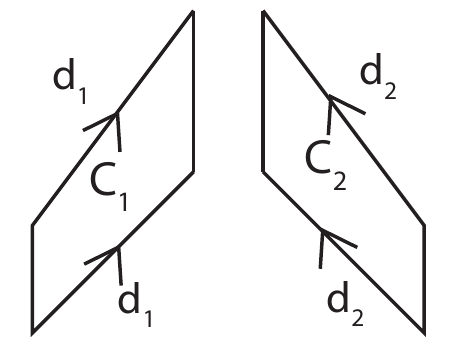}} \caption{The faces
$C_1$ and $C_2$ and the edges $d_1$ and $d_2$. } \label{SlrCZoomup}
\end{center}
\end{figure}

\begin{figure}
\begin{center}
\scalebox{.8}{\includegraphics{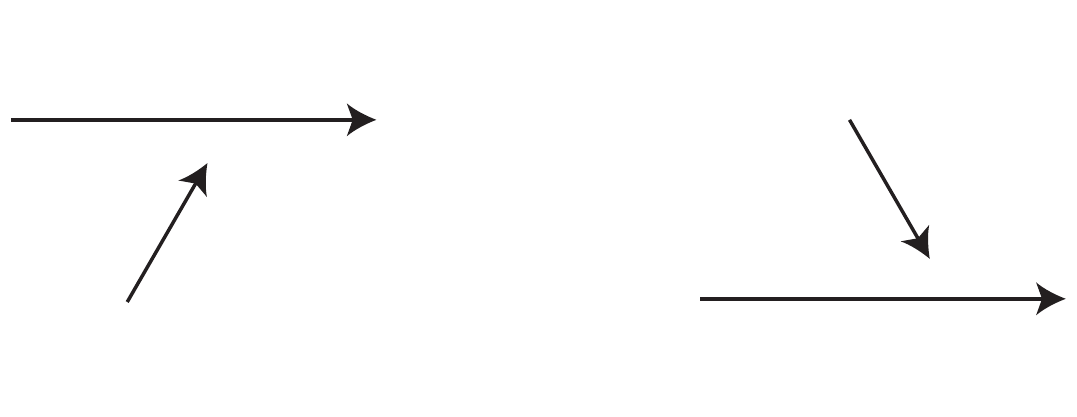}} \caption[Angles of the vector field.]{The angles
that the vector field makes with the edges $d_1$ and $d_2$ are
$\frac{4 \pi}{3}$ and $\frac{2 \pi}{3}$. } \label{SlrCAngles}
\end{center}
\end{figure}

Faces $C_1$ and $C_2$ are slightly different, having opposite
orientation.  See Figure \ref{SlrCZoomup}.  The edge $d_1$ gets sent
to a line of depth $\frac{4 \pi}{3}$, while $d_2$ gets sent to a
line at depth $\frac{2\pi}{3}$.  This is shown in Figure
\ref{SlrCAngles}. The faces will be glued together in a picture identical to Figure \ref{SlrBGluing}, with the appropriate labels (and the faces separated). A careful check will reveal that the map on the common edge of $C_1$ and $B_1$ is well defined, as are all maps on `vertical' edges in the gluing. This is important, as many maps can be created that seem to give the unit tangent bundle, but aren't well defined on the vertical edges.

\begin{figure}
\begin{center}
\scalebox{.8}{\includegraphics{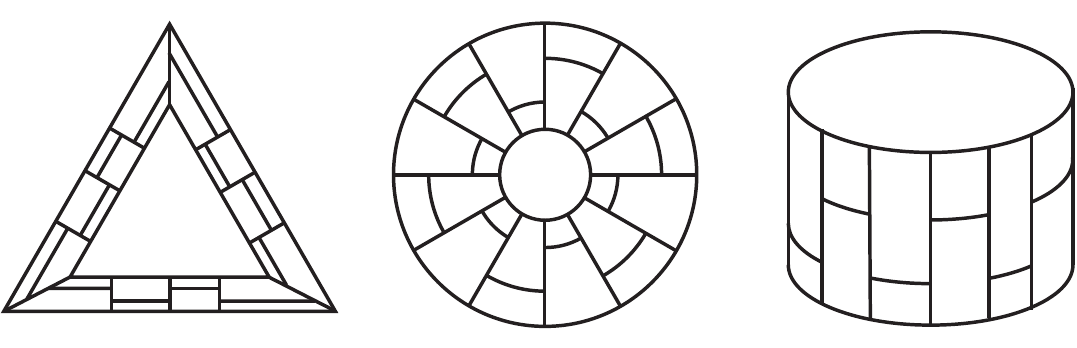}} \caption[The gluing map.]{The gluing
map for M. The first two are just homotopies of each other, and the
last is a 3-d figure of the fundamental polyhedron.}
\label{SlrGluing}
\end{center}
\end{figure}

\begin{figure}
\begin{center}
\scalebox{.8}{\includegraphics{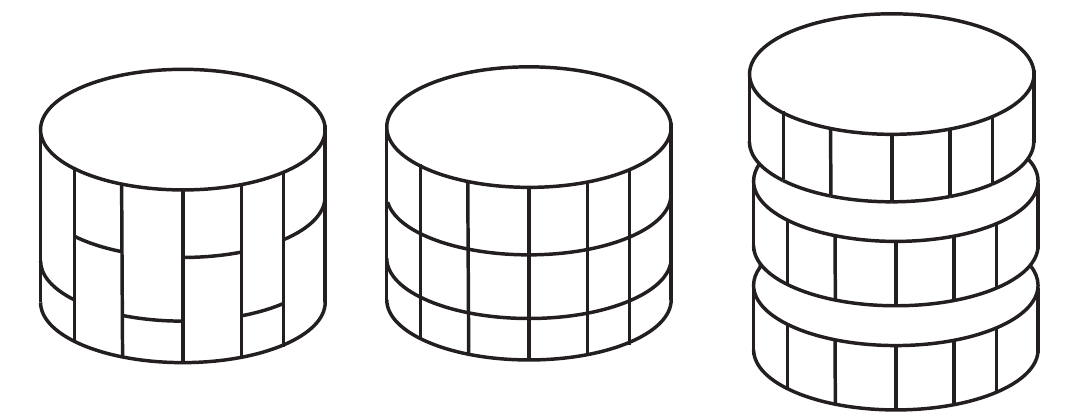}} \caption[Slicing the manifold.]{We slice
the fundamental polyhedron into three smaller cylinders. This makes
each horizontal edge have valence four and each vertical edge have
valence three.} \label{SlrSlicing}
\end{center}
\end{figure}

Continuing all along, we get the gluing map shown in Figure
\ref{SlrGluing}. This can be simplified further by taking the
resulting `cylinder' and slicing horizontally, adding new lines. See
Figure \ref{SlrSlicing}.  Each `vertical' edge has edge cycle length 3, and
each horizontal edge has edge cycle length 4.

Note that the three resulting polyhedra are combinatorially identical to the polyhedron we
used for $N\times S^1$ in Section \ref{H2RSection}.  However,
in $N\times S^1$, each polyhedron corresponded to a single
group element, while in our manifold M, a group element corresponds
to three polyhedra.  Also, in the universal cover of M, neighboring
group elements in the `horizontal' direction are not at the same height; the middle layer of one group
of three polyhedra will get glued alternately to the top and bottom
layers of the neighboring groups of three polyhedra.

The reason that manifolds from these two geometries can have the same subdivision rules is that the two geometries are quasi-isometric.

\section{\texorpdfstring{$\mathbb{S}^3$ geometry: The 3-sphere}{S\textthreesuperior geometry: The 3-sphere}}

Each manifold in this geometry is finitely covered by $\mathbb{S}^3$, which has empty boundary. Thus, these manifolds all have an `empty' subdivision rule. The space at infinity is the empty set, which is never subdivided.

This is different from having no subdivision rule. Groups in this geometry are Gromov hyperbolic (because they are finite), and their hyperbolic boundary is the empty set. Our subdivision rule is the identity on the empty set.

\section{Sol geometry: Nonexistence of subdivision rules}
There are no known subdivision rules for the last two geometries, Nil and Sol.

\begin{figure}
\begin{center}
\scalebox{.4}{\includegraphics{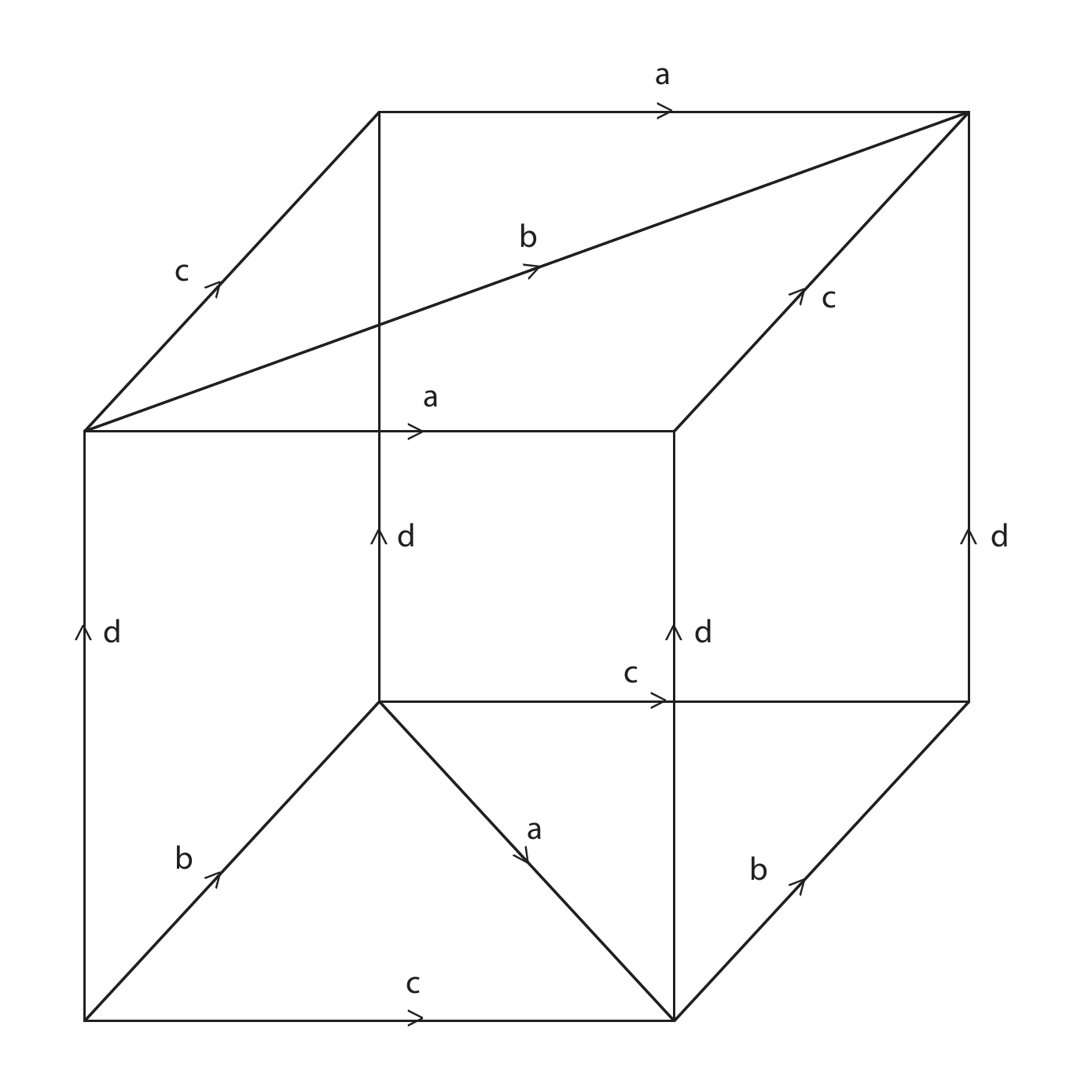}} \caption{A fundamental domain for a solv manifold} \label{SolvCube}
\end{center}
\end{figure}

The standard manifolds for Sol groups are torus bundles over the circle with hyperbolic gluing map. see Figures \ref{NilCube} and \ref{SolvCube}.

The existence of subdivision rules for Solv manifolds is related to the notion of almost convexity. This definition follows \cite{AlmCon}. Here, $S(n)$ and $B(n)$ refer to spheres and balls of radius $n$ in the Cayley graph under the given metric.

\begin{defi}
Let $d(x,y)$ be the word metric in a Cayley graph of a group $G$ with a given generating set. The group $G$ is \textbf{almost convex} ($m$) with respect to this generating set if there is a constant $K(m)$ with the following property such that any two elements $g, g'$ in $S(n)$ with $d(g,g')\leq m$ are connected by an edge path in $B(n)$ of length bounded by $K(m)$. We say that $G$ is \textbf{almost convex} with respect to this generating set if it is almost convex ($m$) for all $m$. A group is almost convex if it is almost convex with regards to some generating set.
\end{defi}

Cannon has shown \cite{AlmCon} that a group with a given generating set is almost convex if it is almost convex (2) with that generating set.

Almost convex groups are exactly the groups which can be constructed efficiently by a local replacement rule on the Cayley graph \cite{ConMod}. A local replacement rule in this sense is more general than our finite replacement rule of Definition \ref{ReplaceDef}, because it deals only with the Cayley graph itself, while ours requires the cell structure of each $S(n)$ to be that of 2-dimensional sphere. However, our finite replacement rules are replacement rules in the sense of \cite{ConMod}. Thus, only almost convex groups can have a finite replacement rule.

Cannon, Floyd, and Parry have shown that Sol groups cannot be almost convex with respect to any generating set \cite{Solv}. This shows that Sol manifolds do not have a finite replacement rule that is associated to the Cayley graph in the way that the other rules we obtained are.

Thus, we cannot obtain subdivision rules for these manifolds by first creating a replacement rule that constructs the universal cover, as we have done. While it may be possible to create some form of subdivision rule associated to these manifolds, it could not be directly connected to the Cayley graph as our other subdivision rules are.

This concludes the proof of Theorem \ref{BigTheorem}.

\section{Nil geometry}\label{NilSection}
\begin{figure}
\begin{center}
\scalebox{.4}{\includegraphics{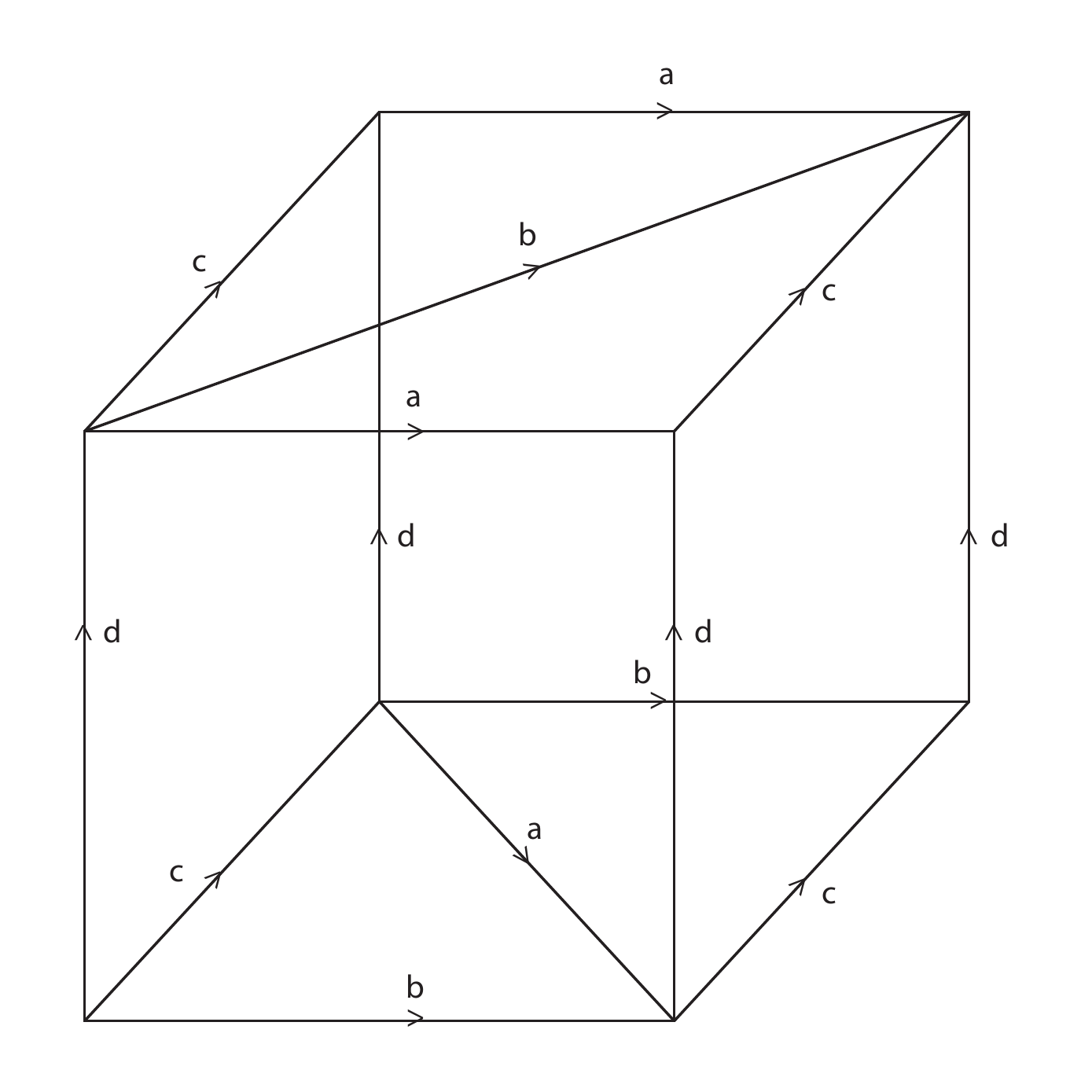}} \caption{A fundamental domain for a nil manifold} \label{NilCube}
\end{center}
\end{figure}

The standard manifold for Nil geometry is the quotient of the 3-dimensional real Heisenberg group by the integral subgroup. As for Sol geometry, it is easy to find a fundamental domain for this manifold (see Figure \ref{NilCube}). Unlike Sol geometry, Nil manifold groups are almost convex \cite{Nil}, meaning that their Cayley graphs with some generating set can be constructed by a local replacement rule in some sense.

However, we again find difficulty creating a finite replacement rule (recall Definition \ref{ReplaceDef}), because the standard generating sets for the standard Nil manifold all have infinitely many cone types \cite{Cone}.

\begin{defi} The \textbf{shadow} of a vertex $v$ in a Cayley graph is the set of all vertices and edges which can be reached by a geodesic ray from the identity going through $v$. Two vertices in a Cayley graph are said to have the same \textbf{cone type} if their shadows are isomorphic.
\end{defi}

In all of our examples of finite replacement rules, the tiles in the tiling of the sphere corresponded exactly to group elements, and the Cayley graph could be reconstructed from the finite replacement rule alone by taking the dual graph of each stage of the replacement rule and adding vertical lines from the vertex associated to a tile to the vertex associated to all tiles that it is replaced by. It is easy to see that cone types in the Cayley graph can be put in one-to-one correspondence with tile types of the finite replacement rule. Thus, the standard Nil manifold with the standard generating sets cannot have a finite replacement rule associated to its Cayley graph, and thus we cannot obtain a subdivision rule for this manifold in the same way that we did above.

We can begin constructing the universal cover, but groups of tiles collapse together in a growing variety of ways.

It is possible that there is some Nil manifold group with some other generating set that has finitely many cone types, but it seems unlikely. Thus, the existence of Nil manifold subdivision rules is unresolved.
\section{Summary}\label{summary}

We summarize the existence of subdivision rules for the eight geometries. We list them by how close their manifolds are to having a conformal subdivision rule on the sphere with combinatorial mesh going to 0 (the original, hyperbolic setting).

Sol: No sol groups are almost convex. No finite subdivision rules coming from finite replacement rules associated to the Cayley graph can exist.

Nil: Nil groups are almost convex, but the standard generating sets for the integral Heisenberg group have infinitely many cone types. No finite subdivision rules coming from finite replacement rules associated to the Cayley graph can exist.

$\mathbb{S}^3$: The model manifold $\mathbb{S}^3$ has empty boundary with an empty subdivision rule.

$\mathbb{S}^2 \times \mathbb{R}$: The model manifold $\mathbb{S}^2 \times \mathbb{S}^1$ has a constant subdivision rule acting on two boundary spheres.

$\mathbb{E}^3$: The model manifold $\mathbb{S}^1 \times \mathbb{S}^1 \times \mathbb{S}^1$ has a subdivision rule on the sphere with combinatorial mesh approaching 0 only at finitely many points.

$\mathbb{H}^2 \times \mathbb{R}$: The model manifold $N_{-1} \times \mathbb{S}^1$ has a subdivision rule on the sphere with combinatorial mesh going to 0 only along a 1-dimensional subset.

$\widetilde{SL_2}(\R)$: Same as $\mathbb{H}^2 \times \mathbb{R}$, except the model manifold is the unit tangent bundle of $N_{-1}$.

$\mathbb{H}^3$: These manifolds have conformal subdivision rules on the sphere with combinatorial mesh going to 0 everywhere.

This is summarized in Figure \ref{AllSubs}.

\section{Future work}
As stated before, the subdivision rules we have developed are highly dependent on the chosen manifold and the chosen generating set. In future work, we hope to show that all subdivision rules for manifolds in a given geometry share certain characteristics.

Two questions we would like to answer in the future are:

Do any generating sets for Nil manifolds have a subdivision rule that constructs the universal cover?

Is there a subdivision rule different from the type considered above that represents a Sol manifold in some sense?

\bibliographystyle{plain}
\bibliography{ManifoldPaper}

\end{document}